\documentclass[10pt]{amsart}
\usepackage{indentfirst,latexsym,bm}
\usepackage{amsfonts}
\usepackage{amssymb}
\usepackage{amsmath}
\usepackage{amsbsy}
\usepackage{dsfont}
\usepackage{amsthm}
\usepackage{hyperref}
\usepackage{amscd}
\usepackage[all]{xy}
\usepackage{chemarrow}
\usepackage{multirow}
\usepackage[titletoc]{appendix}
\usepackage{enumitem}
\begin{document}
%\begin{CJK}{GBK}{song}
\renewcommand{\baselinestretch}{1.2}
\thispagestyle{empty}
%\vspace*{2in}
\title[Hopf algebras over basic Hopf algebras]
{Finite-dimensional Hopf algebras over the smallest non-pointed basic Hopf algebra}
\author{Rongchuan Xiong}
\address{School of Mathematical Sciences, Shanghai Key Laboratory of PMMP,East China Normal University, Shanghai 200241, China}
\email{rcxiong@foxmail.com}
%\thanks{}
%\subjclass[2010]{16T05, 16S35, 18D10}
\date{\today}
\maketitle

\newtheorem{question}{Question}
\newtheorem{defi}{Definition}[section]
\newtheorem{conj}{Conjecture}
\newtheorem{thm}[defi]{Theorem}
\newtheorem{lem}[defi]{Lemma}
\newtheorem{pro}[defi]{Proposition}
\newtheorem{cor}[defi]{Corollary}
\newtheorem{rmk}[defi]{Remark}
\newtheorem{Example}{Example}[section]

\theoremstyle{plain}
\newcounter{maint}
\renewcommand{\themaint}{\Alph{maint}}
\newtheorem{mainthm}[maint]{Theorem}

\theoremstyle{plain}
\newtheorem*{proofthma}{Proof of Theorem A}
\newtheorem*{proofthmb}{Proof of Theorem B}
\newcommand{\tabincell}[2]{\begin{tabular}{@{}#1@{}}#2\end{tabular}}

\newcommand{\C}{\mathcal{C}}
\newcommand{\D}{\mathcal{D}}
\newcommand{\A}{\mathcal{A}}
\newcommand{\De}{\Delta}
\newcommand{\M}{\mathcal{M}}
\newcommand{\K}{\mathds{k}}
\newcommand{\E}{\mathcal{E}}
\newcommand{\Pp}{\mathcal{P}}
\newcommand{\Lam}{\lambda}
\newcommand{\As}{^{\ast}}
\newcommand{\Aa}{a^{\ast}}
\newcommand{\Ab}{(a^2)^{\ast}}
\newcommand{\Ac}{(a^3)^{\ast}}
\newcommand{\Ad}{(a^4)^{\ast}}
\newcommand{\Ae}{(a^5)^{\ast}}
\newcommand{\B}{b^{\ast}}
\newcommand{\BAa}{(ba)^{\ast}}
\newcommand{\BAb}{(ba^2)^{\ast}}
\newcommand{\BAc}{(ba^3)^{\ast}}
\newcommand{\BAd}{(ba^4)^{\ast}}
\newcommand{\BAe}{(ba^5)^{\ast}}
\newcommand{\cF}{\mathcal{F}}

\newcommand{\CYD}{{}^{\C}_{\C}\mathcal{YD}}
\newcommand{\AYD}{{}^{\A}_{\A}\mathcal{YD}}
\newcommand{\HYD}{{}^{H}_{H}\mathcal{YD}}
\newcommand{\ydH}{{}^{H}_{H}\mathcal{YD}}
\newcommand{\KYD}{{}^{\mathcal{K}}_{\mathcal{K}}\mathcal{YD}}
\newcommand{\DM}{{}_{D}\mathcal{M}}
\newcommand{\BN}{\mathcal{B}}
\newcommand{\NA}{\mathcal{B}}

\newcommand{\Ga}{g^{\ast}}
\newcommand{\Gb}{(g^2)^{\ast}}
\newcommand{\Gc}{(g^3)^{\ast}}
\newcommand{\Gd}{(g^4)^{\ast}}
\newcommand{\Ge}{(g^5)^{\ast}}
\newcommand{\X}{x^{\ast}}
\newcommand{\GXa}{(gx)^{\ast}}
\newcommand{\GXb}{(g^2x)^{\ast}}
\newcommand{\GXc}{(g^3x)^{\ast}}
\newcommand{\GXd}{(g^4x)^{\ast}}
\newcommand{\GXe}{(g^5x)^{\ast}}
\newcommand{\I}{\mathds{I}}
\newcommand{\Z}{\mathds{Z}}
\newcommand{\N}{\mathds{N}}
\newcommand{\G}{\mathcal{G}}
\newcommand{\roots }{\boldsymbol{\Delta }}

\newcommand\ad{\operatorname{ad}}
\newcommand{\Alg}{\Hom_{\text{alg}}}
\newcommand\Aut{\operatorname{Aut}}
\newcommand{\AuH}{\Aut_{\text{Hopf}}}
\newcommand\coker{\operatorname{coker}}
\newcommand\car{\operatorname{char}}
\newcommand\Der{\operatorname{Der}}
\newcommand\diag{\operatorname{diag}}
\newcommand\End{\operatorname{End}}
\newcommand\id{\operatorname{id}}
\newcommand\gr{\operatorname{gr}}
\newcommand\GK{\operatorname{GKdim}}
\newcommand{\Hom}{\operatorname{Hom}}
\newcommand\ord{\operatorname{ord}}
\newcommand\rk{\operatorname{rk}}
\newcommand\soc{\operatorname{soc}}

\newcommand{\bp}{\mathbf{p}}
\newcommand{\bq}{\mathbf{q}}
\newcommand\Sb{\mathbb S}
\newcommand\cR{\mathcal{R}}
\newcommand{\grAYD}{{}^{\gr\A}_{\gr\A}\mathcal{YD}}
\newcommand{\Dchaintwo}[3]{\xymatrix@C-4pt{\overset{#1}{\underset{x }{\circ}}\ar
@{-}[r]^{#2}
& \overset{#3}{\underset{v_1 }{\circ}}}}
% p1  p2  p3
% x-------x

\newcommand{\Dchaintwoa}[3]{\xymatrix@C-4pt{\overset{#1}{\underset{\  }{\circ}}\ar
@{-}[r]^{#2}
& \overset{#3}{\underset{\ }{\circ}}}}
% p1  p2  p3
% x-------x

\newcommand{\Dchainthree}[8]{\xymatrix@C-2pt{
\overset{#1}{\underset{#2}{\circ}}\ar  @ {-}[r]^{#3}  & \overset{#4}{\underset{#5
}{\circ}}\ar  @{-}[r]^{#6}
& \overset{#7}{\underset{#8}{\circ}} }}
% p1  p3  p4  p6  p7
% x-------x-------x
% p2      p5      p8

\newcommand{\Dtriangle}[6]{
\xymatrix@R-12pt{  &    \overset{#1}{\underset{x}{\circ}} \ar  @{-}[dl]_{#2}\ar  @{-}[dr]^{#3} & \\
\overset{#4}{\underset{e_1}{\circ}} \ar  @{-}[rr]^{#5}  &  &\overset{#6}{\underset{v_1}{\circ}} }}
%         p1
%         x
%     p2 / \ p3
%       /   \
%   p4 x-----x p6
%         p5

\begin{abstract}
We classify finite-dimensional Hopf algebras over an algebraically closed field of characteristic zero whose Hopf coradcial is isomorphic to the smallest non-pointed basic Hopf algebra, under the assumption that the diagrams are strictly graded.  In particular, we obtain some new Nichols algebras of non-diagonal type and new finite-dimensional Hopf algebras without the dual Chevalley property.

\end{abstract}

\section{Introduction}
Let $\K$ be an algebraically closed field of characteristic zero. This work is a continuation of the paper \cite{GG16} on the classification of finite-dimensional Hopf algebras over $\K$ without the dual Chevalley property, that is, the coradical is not a subalgebra. Until now, there are few concrete examples of such Hopf algebras without pointed duals, with some exceptions in \cite{GG16,HX17,X17}. Let $\mathcal{K}$ be the smallest  Hopf algebra without the dual Chevalley property. It is basic with the dual a Radford algebra $\A$ $($see \eqref{eqDefA} for the definition$)$.  The authors in \cite{GG16} determined all finite-dimensional Hopf algebras over $\mathcal{K}$ whose diagrams are Nichols algebras over the indecomposable objects in $\KYD$ via the generalized lifting method proposed by Andruskiewitsch and Cuadra \cite{AC13}.

The generalized lifting method is a generalization of the \emph{lifting method }introduced by Andruskiewitsch-Schneider \cite{AS98}. Let $A$ be a Hopf algebra without the dual Chevalley property. Andruskiewitsch-Cuadra replaced the coradical filtration  $\{A_{(n)}\}_{n\geq 0}$ with the standard filtration $\{A_{[n]}\}_{n\geq 0}$, which is
defined recursively by $A_{[n]}=A_{[n-1]}\bigwedge A_{[0]}$, where $A_{[0]}$ called the Hopf coradical of $A$ is the subalgebra generated by the coradical $A_0$.
Assume that $S_A(A_{[0]})\subseteq A_{[0]}$, it turns out that
the associated graded coalgebra
$\gr A=\oplus_{n=0}^{\infty}A_{[n]}/A_{[n-1]}$ with $A_{[-1]}=0$
is a Hopf algebra and $\gr A\cong R\sharp A_{[0]}$ as Hopf algebras, where $R=(\gr A)^{co A_{[0]}}=\oplus_{n\geq 0}R(n)$ called the \emph{diagram} of $A$ is a connected $\mathds{N}$-graded braided Hopf algebra
in ${}^{A_{[0]}}_{A_{[0]}}\mathcal{YD}$. Moreover, $R(1)$ called the \emph{infinitesimal braiding} of $A$ is a subspace of $\Pp(R)$. If the coradical $A_0$ is a Hopf subalgebra, then the standard filtration coincides with the coradical
one.  In this case, $\gr A$ is coradically graded and the diagram $R$ of $A$ is \emph{strictly graded}, that is,
$R(0)=\K,R(1)=\Pp(R)$. In general, it is open whether the diagram $R$ is strictly graded or not. See \cite{AS02,AC13} for details. To construct Hopf algebras via the generalized lifting method, the following questions are considered (see \cite{AC13}):
\begin{itemize}
  \item {\text{\rm Question\,\uppercase\expandafter{\romannumeral1}. }} Let $C$ be a cosemisimple coalgebra and $\mathcal{S} : C \rightarrow C$
  an injective anti-coalgebra morphism. Classify all Hopf algebras $L$ generated by $C$, such that $S|_C = \mathcal{S}$.
  \item {\text{\rm Question\,\uppercase\expandafter{\romannumeral2}. }} Given $L$ as in the previous item, classify all connected graded
  Hopf algebras $R$ in ${}_L^L\mathcal{YD}$.
  \item {\text{\rm Question\,\uppercase\expandafter{\romannumeral3}. }} Given $L$ and $R$ as in previous items, classify all liftings,
  that is, classify all Hopf algebras $A$ such that $\gr A\cong R\sharp L$. We call $A$ a lifting of $R$ over $L$.
\end{itemize}

 As the aforementioned, the authors in \cite{GG16} determined all finite-dimensional Hopf algebras such that the diagrams $R$ are Nichols algebras $\BN(V)$ over the indecomposable objects $V$ in the case $L=\mathcal{K}$. We continue to study  these questions \uppercase\expandafter{\romannumeral2} \& \uppercase\expandafter{\romannumeral3} in the case that $V$ are semisimple objects in $\KYD$.

The Hopf algebra $\mathcal{K}$ is defined in Proposition \ref{proStrucOfC}. Notice it that the description of $\mathcal{K}$ is different from that in \cite[Proposition\,2.1.]{GG16}. We describe the structure and representation of the Drinfeld double $\D:=\D(\mathcal{K}^{cop})$, and describe the simple objects  in $\KYD$ by using the equivalence $\KYD\cong {}_{\D(\mathcal{K}^{cop})}\mathcal{M}$
\cite[Proposition\,10.6.16]{M93}.  Indeed, as stated in \cite{GG16},  there are $4$ one-dimensional objects
$\K_{\chi_{k}}$ with $0\leq k<4$ and $12$ two-dimensional objects $V_{i,j}$ with
$(i,j)\in\Lambda=\{(i,j)\in\N^2\mid 0\leq i,j<4, 2i\neq j \text{~mod~}4\} $. Then we determine all finite-dimensional Nichols algebras $\BN(V)$ in $\KYD$ and  show that they can be obtained by splitting Nichols algebras of diagonal type $($standard type A$)$. The main result is given as follows.
\begin{mainthm}\label{thmA}
Let $V\in\KYD$ such that the Nichols algebra $\BN(V)$ is finite-dimensional, then $V$ is isomorphic either to $\oplus_{k=1}^n\K_{\chi^{i_k}}$ with $i_k\in\{1,3\}$, $V_{1,j}$, $V_{2,j}$, $V_{1,j}\oplus \K_{\chi}$, $V_{2,j}\oplus\K_{\chi^3}$, $V_{1,1}\oplus V_{1,3}$ or $V_{2,1}\oplus V_{2,3}$, where $n<\infty$ and $j\in\{1,3\}$. Moreover, the generators and relations of $\BN(V)$ are given as follows:
\center
%\caption{Finite-dimensional Nichols algebras}
\begin{tabular}{|c|c|c|c|}
  \hline
   $V$           &   relations of $\BN(V)$ with generators $x,y,z,t$         & $\dim\BN(V)$  \\\hline
$\oplus_{k=1}^n\K_{\chi^{i_k}}$   &   $\otimes_{1\leq k\leq n}\bigwedge \K_{\chi^{i_k}}$  &   $2^{n}$   \\\hline
$V_{1,j}$       &   $x^2+2\xi y^2=0, xy+yx=0, x^4=0  $    &  8  \\\hline
$V_{2,j}$       &    $ x^2=0, yx+\xi^jxy=0, y^4=0$   &  8  \\\hline
$V_{1,j}\oplus \K_{\chi}$ &\tabincell{c}{ $x^4=0$, $xy+yx=0$, $x^2+2\xi y^2=0$, $z^2=0$,$zx^2+(1-\xi^j)xzx-\xi^j x^2z=0$, \\$\xi^j xyz-\xi^j xzy+yzx+zxy=0$,
$\frac{1}{2}\xi(1+\xi^{-j})(xz)^2(yz)^2+(yz)^4+(zy)^4=0$. }& 128  \\\hline
$V_{2,j}\oplus\K_{\chi^{3}}$&\tabincell{c}{ $x^2=0$, $yx+\xi^j xy=0$, $y^4=0$, $z^2=0$,\\$xyz+xzy+yzx+zyx=0$, $(zx)^4+(xz)^4=0$,\\$\frac{1}{2}\xi(1+\xi^{-j})xzx-\xi^j y^2z+(1-\xi^j)yzy+zy^2=0$. } & 128        \\\hline
$V_{1,1}\oplus V_{1,3}$&\tabincell{c}{ $x^4=0$, $xy+yx=0$, $x^2+2\xi y^2=0$, $zx-\xi xz=0$,\\ $tx+zy+yz+xt=0$,
 $ty+yt+\frac{1}{2}(1-\xi)xz=0$,\\ $z^4=0$, $zt+tz=0$, $z^2+2\xi t^2=0$} &       128              \\\hline
 $V_{2,1}\oplus V_{2,3}$&\tabincell{c}{ $x^2=0$, $yx+\xi xy=0$, $y^4=0$,  $zx+xz=0$,\\ $tx-zy-\xi yz+\xi xt=0$,
 $ty-\xi yt+\frac{1}{2}(\xi-1)xz=0$,\\ $z^2=0$, $tz-\xi zt=0$, $t^4=0$} &       128            \\\hline
\end{tabular}
%\begin{tabular}{|c|c|c|c|}
%  \hline
%   $V$           &   relations of $\BN(V)$ with generators $x,y,z,t$         & $\dim\BN(V)$ & $\BN(U)$ \\\hline
%$\oplus_{k=1}^n\K_{\chi^{i_k}}$   &   $\otimes_{1\leq k\leq n}\bigwedge \K_{\chi^{i_k}}$  &   $2^{n}$ &$A_1\times\cdots\times A_1$  \\\hline
%$V_{1,j}$       &   $x^2+2\xi y^2=0, xy+yx=0, x^4=0  $    &  8  & standard type $A_2$\\\hline
%$V_{2,j}$       &    $ x^2=0, yx+\xi^jxy=0, y^4=0$   &  8  & standard type $A_2$ \\\hline
%$V_{1,j}\oplus \K_{\chi}$ &\tabincell{c}{ $x^4=0$, $xy+yx=0$, $x^2+2\xi y^2=0$, $z^2=0$,\\$zx^2+(1-\xi^j)xzx-\xi^j x^2z=0$, \\$\xi^j xyz-\xi^j xzy+yzx+zxy=0$,\\
%$\frac{1}{2}\xi(1+\xi^{-j})(xz)^2(yz)^2+(yz)^4+(zy)^4=0$. }& 128  & standard type $A_3$\\\hline
%$V_{2,j}\oplus\K_{\chi^{3}}$&\tabincell{c}{ $x^2=0$, $yx+\xi^j xy=0$, $y^4=0$, $z^2=0$,\\$xyz+xzy+yzx+zyx=0$, $(zx)^4+(xz)^4=0$,\\$\frac{1}{2}\xi(1+\xi^{-j})xzx-\xi^j y^2z+(1-\xi^j)yzy+zy^2=0$. } & 128    & standard type $A_3$      \\\hline
%$V_{1,1}\oplus V_{1,3}$&\tabincell{c}{ $x^4=0$, $xy+yx=0$, $x^2+2\xi y^2=0$, $zx-\xi xz=0$,\\ $tx+zy+yz+xt=0$,
% $ty+yt+\frac{1}{2}(1-\xi)xz=0$,\\ $z^4=0$, $zt+tz=0$, $z^2+2\xi t^2=0$} &       128    & standard type $A_3$          \\\hline
% $V_{2,1}\oplus V_{2,3}$&\tabincell{c}{ $x^2=0$, $yx+\xi xy=0$, $y^4=0$,  $zx+xz=0$,\\ $tx-zy-\xi yz+\xi xt=0$,
% $ty-\xi yt+\frac{1}{2}(\xi-1)xz=0$,\\ $z^2=0$, $tz-\xi zt=0$, $t^4=0$} &       128 &  standard type $A_3$            \\\hline
%\end{tabular}
\end{mainthm}
The Nichols algebras $\BN(V_{i,j})$ have  appeared in \cite{GG16} and were shown that they are not of diagonal type. The Nichols algebras over the semisimple objects are also of non-diagonal type except for those over the direct sum of one-dimensional objects. To the best of our knowledge, they constitute new examples of finite-dimensional Nichols algebras. Moreover, their generalized Cartan matrices are of type $A_2$ and they can give rise to a Weyl groupoid of standard type $A_2$ $($see Corollary \ref{corCartan}$)$.

The proof of Throrem \ref{thmA} uses the equivalence $\KYD\cong\AYD$ by \cite[Proposition\,2.2.1.]{AG99} and depends highly on the classification results of finite-dimensional Nichols algebras of diagonal type $($see \cite{H09} and the references therein$)$. As $\KYD\cong\AYD$,  the braided vector space $V$ in $\KYD$ can be regarded as the corresponding object in $\AYD$ if there is no confusion. Notice that $\BN(V)\sharp\A$ is a pointed Hopf algebra. In order to show that $\dim\BN(V)=\infty$ for some $V\in\KYD$, we turn to find an Hopf subalgebra  $B$ of $\BN(V)\sharp\A$ and show that  $\gr B$ is infinite-dimensional by using the results in \cite{H09}.  Hence we can discard Nichols algebras of infinite dimension, then prove the remaining are finite-dimensional by computing their relations and PBW bases.
See section \ref{secNicholsalgebra} for more details. Here we also refer to a recent work \cite{AA} for a characterization of finite-dimensional Nichols algebras over basic Hopf algebras.

Finally, we study the liftings of the Nichols algebras in Theorem \ref{thmA} following the techniques in \cite{AS98,GG16}. It turns out that the Nichols algebra $\BN(V)$  does not admit non-trivial deformations, where $V$ is isomorphic either to $\oplus_{i=1}^n\K_{\chi^{n_i}}$ with $n_i\in\{1,3\}$, $V_{2,1}$, $V_{2,3}$, $V_{2,1}\oplus\K_{\chi^3}$, $V_{2,3}\oplus\K_{\chi^3}$ or $V_{2,1}\oplus V_{2,3}$. The bosonizations of them are basic. The remaining admit non-trivial deformations. Hence we define five families of Hopf algebras $\mathfrak{A}_{1,j}(\mu)$, $\mathfrak{A}_{1,j,1}(\mu,\nu)$ for $j=1,3$ and $\mathfrak{A}_{1,1,1,3}(\mu,\nu)$ $($see Definitions \ref{defiV13} \ref{defiV131} \& \ref{defiV1313}$)$, and show that they are liftings of the Nichols algebras $\BN(V_{1,j})$, $\BN(V_{1,j}\oplus\K_{\chi})$ for $j=1,3$ and $\BN(V_{1,1}\oplus V_{1,3})$, respectively. Moreover, under the assumption that the diagrams are strictly graded, we show that all finite-dimensional Hopf algebras over $\mathcal{K}$ are generated in degree one with respect to the standard filtration $($see Theorem \ref{thmGeneratedbyone}$)$. As a summary, we have the following
\begin{mainthm}\label{thmB}
Let $A$ be a finite-dimensional Hopf algebra over $\mathcal{K}$. Assume that the diagram of $A$ is strictly graded. Then $A$ is isomorphic either to $\BN(\oplus_{k=1}^n\K_{\chi^{i_k}})\sharp\mathcal{K}$, $\mathfrak{A}_{1,j}(\mu)$, $\BN(V_{2,j})\sharp\mathcal{K}$, $\mathfrak{A}_{1,j,1}(\mu,\nu)$, $\BN(V_{2,j}\oplus\K_{\chi^3})\sharp\mathcal{K}$, $\mathfrak{A}_{1,1,1,3}(\mu,\nu)$ or $\BN(V_{2,1}\oplus V_{2,3})\sharp\mathcal{K}$ for $j\in\{1,3\},\mu,\nu\in\K$.
\end{mainthm}

The Hopf algebras $\mathfrak{A}_{1,j}(\mu)$ with $j=1,3$ have already appeared in \cite{GG16} up to isomorphism. They are $64$-dimensional Hopf algebras without the dual Chevalley property. It should be figured out that the structure of $\mathfrak{A}_{1,j,1}(\mu,\nu)$ is not completely determined if $\nu\neq 0$. Due to complicated commutators relations and computations, the lifting of the relation $\frac{1}{2}\xi(1+\xi^{-j})(xz)^2(yz)^2+(yz)^4+(zy)^4=0$ of  $\BN(V_{1,j}\oplus\K_{\chi})$ is not clear in general but characterized by the parameters $\mu,\nu$. $\mathfrak{A}_{1,j,1}(\mu,0)$ for $j=1,3$ and $\mathfrak{A}_{1,1,1,3}(\mu,\nu)$ also do not have the dual Chevalley property. They are not basic and constitute new examples of finite-dimensional Hopf algebras except for $\mu=0=\nu$. See section $\ref{secHopfalgebra}$ for more details.

The structure of the paper is given as follows. In section $\ref{Preliminary}$ we first introduce some basic definitions and facts about Yetter-Drinfeld modules,
 Nichols algebras, and redescribe the structure of $\mathcal{K}$ and the representation of the Drinfeld double $\D(\mathcal{K}^{cop})$.
In section $\ref{secNicholsalgebra}$, we determine all finite-dimensional Nichols algebras in $\KYD$ and present explicitly them by generators and relations. In section $\ref{secHopfalgebra}$, we mainly compute the liftings of the Nichols algebras in Theorem \ref{thmA} and prove Theorem \ref{thmB}.

\section{Preliminaries}\label{Preliminary}
\subsection*{Conventions}We work over an algebraically closed field $\K$ of characteristic zero and denote
by $\xi$ a primitive $4$-th root of unity. The references for Hopf algebra theory are \cite{M93,R11}.

The notation for a Hopf algebra $H$ over $\K$ is standard: $\Delta$, $\epsilon$, and $S$ denote the comultiplication, the counit and the antipode.
 We use Sweedler's notation for the comultiplication and coaction. %, for example, for any $h\in H$, $\Delta(h)=h_{(1)}\otimes h_{(2)}$.
 %Denote by $H^{op}$
% the Hopf algebra with the opposite multiplication, by $H^{cop}$ the Hopf algebra with the opposite comultiplication, and by $H^{bop}$
% the Hopf algebra $H^{op\,cop}$.
Denote by $\G(H)$ the set of group-like elements of $H$. For any $g,h\in\G(H)$, $\Pp_{g,h}(H)=\{x\in H\mid \Delta(x)=x\otimes g+h\otimes x\}$. In particular, the linear space $\Pp(H):=\Pp_{1,1}(H)$ is called the set
 of primitive elements. %Denote by $\otimes$ for simplicity the tensor products over $\K$ and by ${}_H\mathcal{M}$ the category of left $H$-modules.

 Given two  $($braided monoidal$)$ categories $\mathfrak{C}$ and $\mathfrak{D}$, denote by $\mathfrak{C}\cong\mathfrak{D}$ the $($braided monoidal$)$ equivalence between $\mathfrak{C}$ and $\mathfrak{D}$.
 %If $V$ is a $\K$-vector space, $v\in V$ and $f\in V\As$, we use either $f(v)$, $\langle f$, $v\rangle$, or $\langle v,f\rangle$ to
% denote the evaluation.
 Given $n\geq 0$, we denote $\Z_n=\Z/n\Z$ and $\I_{0,n}=\{0,1,\ldots,n\}$. In particular, the operations $ij$ and $i\pm j$ are considered modulo $n+1$ for $i,j\in\I_{k,n}$ when not specified.

\subsection{ Yetter-Drinfeld modules and Nichols algebras}
Let $H$ be a Hopf algebra with bijective antipode. A left \emph{Yetter-Drinfeld module} $M$ over $H$ is a left $H$-module $(M,\cdot)$ and a left $H$-comodule $(M,\delta)$ satisfying
$
\delta(h\cdot v)=h_{(1)}v_{(-1)}S(h_{(3)})\otimes h_{(2)}\cdot v_{(0)}
$
for all $v\in V,h\in H$.

Let ${}^{H}_{H}\mathcal{YD}$ be the category of Yetter-Drinfeld modules over $H$. Then ${}^{H}_{H}\mathcal{YD}$ is braided monoidal. For $V,W\in {}^{H}_{H}\mathcal{YD}$, the braiding $c_{V,W}$ is given by
\begin{align}\label{equbraidingYDcat}
c_{V,W}:V\otimes W\mapsto W\otimes V, v\otimes w\mapsto v_{(-1)}\cdot w\otimes v_{(0)},\forall\,v\in V, w\in W.
\end{align}
In particular, $(V,c_{V£¬V})$ is a braided vector space, that is, $c:=c_{V,V}$ is a linear isomorphism satisfying the
braid equation $(c\otimes\text{id})(\text{id}\otimes c)(c\otimes\text{id})=(\text{id}\otimes c)(c\otimes\text{id})(\text{id}\otimes c)$. Moreover, ${}^{H}_{H}\mathcal{YD}$ is rigid. The left dual $V\As$ is defined by
\begin{align*}
\langle h\cdot f,v\rangle=\langle f,S(h)v\rangle,\quad f_{(-1)}\langle f_{(0)},v\rangle=S^{-1}(v_{(-1)})\langle f, v_{(0)}\rangle.
%\langle h\cdot f,v\rangle=\langle f,S^{-1}(h)v\rangle,\quad f_{(-1)}\langle f_{(0)},v\rangle=S(v_{(-1)})\langle f, v_{(0)}\rangle.
\end{align*}

Assume that $H$ is a finite-dimensional Hopf algebra. Then by \cite[Proposition\,2.2.1.]{AG99}, $\HYD\cong{}_{H^{\ast}}^{H^{\ast}}\mathcal{YD}$ as braided monoidal categories  via the functor $(F,\eta)$ defined as follows:
 $F(V)=V$ as a vector space,
\begin{align}
\begin{split}\label{eqVHD}
f\cdot v=f(S(v_{(-1)}))v_{(0)},\quad \delta(v)=\sum_{i}S^{-1}(h^i)\otimes h_i\cdot v,\  \text{and}\\
\eta:F(V)\otimes F(W)\mapsto F(V\otimes W), v\otimes w\mapsto w_{(-1)}\cdot v\otimes w_{(0)}
\end{split}
\end{align}
for every $V,W\in\HYD$, $f\in H^{\ast},v\in V, w\in W$. Here $\{h_i\}$ and $\{h^i\}$ are the dual bases of $H$ and $H\As$.

\begin{defi}\cite[Definition\;2.1]{AS02}
Let $H$ be a Hopf algebra and $V\in{}^{H}_{H}\mathcal{YD}$. A braided graded Hopf algebra $R=\oplus_{n\geq 0} R(n)$ in ${}^{H}_{H}\mathcal{YD}$
is called a Nichols algebra over $V$  if
\begin{align*}
R(0)=\mathds{k}, \quad R(1)=V,\quad
R\;\text{is generated as an algebra by}\;R(1),\quad
\Pp(R)=V.
\end{align*}
\end{defi}
Let $V\in\HYD$, we denote by $\BN(V)$ the Nichols algebra over $V$. $\BN(V)$ is unique up to isomorphism and isomorphic to $T(V)/I(V)$, where $I(V)\subset T(V)$ is the largest $\mathds{N}$-graded ideal and coideal in $\HYD$ such that $I(V)\cap V=0$. Moreover,
$\BN(V)$ as a coalgebra and an algebra depends only on $(V,c)$
 and the ideal $I(V)$ is the kernel of the quantum symmetrizer associated to the braiding $c$.

 Let $(W,c)$ be a braided
 vector subspace of $(V,c)$, that is, $W$ is a vector subspace of $V$ such that $c(W\otimes W)\subset W\otimes W$. Then $\BN(W)$
 is a braided Hopf subalgebra of $\BN(V)$. In particular, $\dim \BN(V)=\infty$ if $\dim \BN(W)=\infty$.
See \cite{AS02} for more details.

Let $(V,c)$ be a $n$-dimensional $($rigid$)$ braided vector
space and $f\in V^{\ast}$. Then the \emph{skew-derivation} $\partial_f\in\text{End}\,T(V)$ is given by $\partial_f(1)=0$, $\partial_f(v)=f(v)$ and
%\begin{gather*}
%\partial_f(1)=0,\quad\quad \partial_f(v)=f(v)
$\partial_f(xy)=\partial_f(x)y+\sum_{i}x_i\partial_{f_i}(y),\;\text{where}\;c^{-1}(f\otimes x)=\sum_ix_i\otimes f_i$
%\end{gather*}
for $v\in V,x,y\in T(V)$. Let $\{v_i\}_{1\leq i\leq n}$ and $\{v^i\}_{1\leq i\leq n}$ be the dual bases of $V$ and $V\As$. We denote $\partial_i:=\partial_{v^i}$ for simplicity. The skew-derivation is very useful
to find the relations of Nichols algebra $\BN(V)$: Let $r\in T^m(V)$, $r\in I(V)$ if and only if $\partial_f(r)=0$ for all $f\in V^{\ast}$ if and only if $\partial_i(r)=0$ for all $1\leq i\leq n$.  See \cite[Theorem\,2.9.]{AHS10} for details.

%We close this subsection by giving the explicit relation between $V$ and $V\As$ in ${}^{H}_{H}\mathcal{YD}$.
%\begin{pro}\cite[Proposition\,3.2.30]{AG99}\label{proNicholsdual}
%Let $V$ be an object in ${}_H^H\mathcal{YD}$. If $\BN(V)$ is finite-dimensional, then $\BN(V\As)\cong \BN(V)^{\ast\,op\,cop}$.
%\end{pro}
\subsection{The Weyl groupoid of a Nichols algebra}
 The Weyl groupoid of a Nichols algebra of diagonal type is first introduced in \cite{H06} $($see also \cite{AHS10,HS10,HS10b}$)$. Now we introduce some notations. The braided vector space $(V,c)$ is of \emph{diagonal type} if there is a linear basis $\{x_1,x_2,\ldots, x_n\}$ such that $c(x_i\otimes x_j)=q_{ij}x_j\otimes x_i$ for some $q_{ij}\in\K$. The matrix $\mathbf{q} = (q_{ij})_{i,j\in \I_{1,n}}$ is called the matrix of the braiding. Let $(\alpha_i)_{i\in \I_{1,n}}$ be the canonical basis of $\Z^n$ and $\chi$ the bicharacter on $\Z^n$ such that
$\chi (\alpha _i,\alpha _j)=q_{ij}$ for all $i,j\in \I $. $\BN(V)$ is $\Z^n $-graded with $\deg x_i=\alpha _i$ for all $i\in \I_{1,n}$
and  there is a totally ordered subset $L \subset \BN(V)$ consisting of $\Z^n$-homoge\-neous elements such that a linear basis of $\BN(V)$ is given by
$
	\{ l_1^{m_1}\cdots l_k^{m_k} \,|\,k\in \N _0,l_1>\cdots >l_k\in L,0 < m_i< N_{l_i}\,\text{for all $i \in \I_{1,k}$}\},
$
where $ N_l= \min\{n\in \N_: (n)_{q_{l, l}}=0\} \in \N \cup \{\infty\}. $

 The \emph{generalized Cartan matrix} $(c_{ij})_{i,j\in\I_{1,n}}$ is given by
\begin{align*}
c_{ii} =2,\quad	c_{ij}=-\min\{m\in\N:(m+1)_{q_{ii}}(1-q_{ii}^m q_{ij}q_{ji})=0\},\quad j\neq i.
\end{align*}
Let $s_i\in GL(\Z^n)$ be given by $s_i (\alpha_j) = \alpha_j - c_{ij}\alpha_i$, $j\in \I_{1,n}$.
The reflection at the vertex $i$ of $\bq$ can be given by the new matrix of braiding $\cR^i(\bq) = (t_{jk}^{(i)})_{j,k\in\I}$, where
\begin{align}
	t_{jk}^{(i)}&:= \chi(s_i(\alpha_j), s_i(\alpha_k)) =  q_{jk}q_{ik}^{-c_{ij}}q_{ji}^{-c_{ik}}q_{ii}^{c_{ij} c_{ik}}, \quad j,k\in\I_{1,n}.
\end{align}

\begin{thm}\cite{H06}\label{thmReflection}
	If $\cR^i(M)$ is the braided vector space corresponding to $\cR^i(\bq)$, then $\dim\BN(\cR^i(V))=\dim\BN(V)$.
\end{thm}

Let $M=\oplus_{i=1}^nM_i$ be a finite-dimensional semisimple Yetter-Drinfeld modules. The Weyl groupoid of $\BN(M)$ is introduced in \cite{AHS10} as a generalization of that in \cite{H06}.
Assume that $\BN(M)$ is finite-dimensional for simplicity. Then the generalized Cartan matrix $(c_{ij})_{i,j\in\I_{1,n}}$ is given by $c_{ii}=2$ and
$
-c_{ij}:=\max\{m\in\N\mid \ad_c(M_i)^m(M_j)\neq 0\}
$ and the $i$-reflection of $M$ is defined by $\cR_i(M) = (V_1,\dots,V_{\theta})$, where
\begin{align*}
V_j &= \begin{cases}
M_i^*, &\text{ if } j=i,\\
(\ad_c M_i)^{-c_{ij}}(M_j), &\text{ if } j \neq i.
\end{cases}
\end{align*}
We refer to \cite[subsection\,3.5.]{AHS10} and also \cite{HS10,HS10b} and the references therein for details.
\subsection{Bosonization and Hopf algebras with a projection}
Let $R$ be a Hopf algebra in ${}^{H}_{H}\mathcal{YD}$ and denote the coproduct
by $\Delta_R(r)=r^{(1)}\otimes r^{(2)}$. The bosonization $R\sharp H$ is defined as follows:
$R\sharp H=R\otimes H$ as a vector space, and the multiplication and comultiplication are given by the smash product and smash-coproduct, respectively:
\begin{align}
(r\#g)(s\#h)=r(g_{(1)}\cdot s)\#g_{(2)}h,\quad
\Delta(r\#g)=r^{(1)}\#(r^{(2)})_{(-1)}g_{(1)}\otimes (r^{(2)})_{(0)}\#g_{(2)}.\label{eqSmash}
\end{align}
Clearly, the map $\iota:H\rightarrow R\sharp H, h\mapsto 1\sharp h,\ \forall h\in H$ is injective and the map
$\pi:R\sharp H\rightarrow H,r\sharp h\mapsto \epsilon_R(r)h,\ \forall r\in R, h\in H$ is surjective such that $\pi\circ\iota=id_H$.
Moreover, $R=(R\sharp H)^{coH}=\{x\in R\sharp H\mid (\text{id}\otimes\pi)\Delta(x)=x\otimes 1\}$.

%Let $R, S$ be the Hopf algebrs in ${}^{H}_{H}\mathcal{YD}$ and $f:R\rightarrow S$ be a bialgebra
%morphism in ${}^{H}_{H}\mathcal{YD}$. The map $f\#id:R\#H\rightarrow S\#H$ is defined by $(f\#id)(r\#h)=f(r)\#h,
%\forall r\in R, h\in H$. In fact, $R\rightarrow R\#H$ and $f\mapsto f\# id$ describes a
%functor from the category of Hopf algebras in ${}^{H}_{H}\mathcal{YD}$ and their morphisms
%to the category of usual Hopf algebras.

Conversely, if $A$ is a Hopf algebra and $\pi:A\rightarrow H$ is a bialgebra map admitting a bialgebra
section $\iota:H\rightarrow A$ such that $\pi\circ\iota=\text{id}_H$, then $A\simeq R\#H$, where $R=A^{coH}$ is a Hopf algebra in ${}^{H}_{H}\mathcal{YD}$. See \cite{R11} for more details. %whose Yetter-Drinfeld module and coalgebra structures are given by:
%\begin{gather}
%h\cdot r=h_{(1)}rS_A(h_{(2)}),\;
%\delta(r)=(\pi\otimes id)\Delta_A(r),\;
%\Delta_R(r)=r_{(1)}(\iota S_H(\pi(r_{(2)})))\otimes r_{(3)},\label{eq-R11}\\
%\epsilon_R=\epsilon_A|_R,\;
%S_R(r)=(\iota\pi(r_{(1)}))S_A(r_{(2)}).
%\end{gather}

\subsection{The Hopf algebra $\mathcal{K}$ and the Drinfeld double $\D(\mathcal{K}^{\text{cop}})$}
Pointed Hopf algebras of dimension $8$ over $\K$ were classified
by \c{S}tefan \cite{S99}. It turns out that these Hopf algebras have pointed duals except for one case
denoted by $\A$ whose Hopf algebra structure is given by
\begin{align}
\A:=\langle g,x\mid g^{4}=1, x^2=1{-}g^2, gx=-xg\;\rangle; \quad \De(g)=g\otimes g,\quad \De(x)=x{\otimes} 1+g{\otimes} x.\label{eqDefA}
\end{align}
Then $\gr\A=\BN(W)\sharp \K[\Gamma]$, where $\Gamma\cong \Z_4$ with the generator $g$ and $W:=\K\{x\}\in{}_{\Gamma}^{\Gamma}\mathcal{YD}$ with the Yetter-Drinfeld module structure given by $g\cdot x=-x$ and $\delta(x)=g\otimes x$.
Moreover, $\A$ is a cocycle deformation of $\gr\A$ \cite{GM07} and whence $\AYD\cong \grAYD$ as braided monoidal categories \cite[Theorem\, 2.7]{MO99}.

Note that $\A\cong (\A_4^{\prime\prime})^{cop}$, where the notation $\A_4^{\prime\prime}$ is introduced in \cite{GV10,GG16}.  The Hopf algebra $\mathcal{K}$ as the dual Hopf algebra of $\A$ is the unique (up to isomorphism) Hopf algebra of dimension $8$ without the dual Chevalley property. Indeed, its coradical $\mathcal{K}_0\simeq \K\oplus \K \oplus C$, where $C$ is a simple coalgebra of dimension $4$.
As a Hopf algebra, the structure of $\mathcal{K}$ is given as follows:
\begin{pro}\label{proStrucOfC}
$\mathcal{K}$ as an algebra is generated by the elements $a$ and $b$ satisfying the relations
\begin{align}
a^{4}=1, \quad b^2=0,\quad ba=\xi ab,\label{eqDef1}
\end{align}
and as a coalgebra is defined by
\begin{align}
 \De(a)=a\otimes a+ \xi^{-1}b\otimes ba^2,\quad \De(b)=b\otimes a^{3}+a\otimes b,\quad \epsilon(a)=1,\quad\epsilon(b)=0,\label{eqDef}
\end{align}
and the antipode is given by $S(a)=a^{-1}$, $ S(b)=\xi^{3}b$.
\end{pro}
\begin{proof}
Similar to the proof of \cite[Lemma 3.3.]{GV10}
\end{proof}

\begin{rmk}\label{rmkDDDDD}
\begin{enumerate}
\item A straight computation shows that $\G(\C)=\{1,a^2\}$, $\Pp_{1,a^2}(\C)=\K\{1-a^2,ba\}$. In particular, the subalgebra generated by $a^2$ and $ba$ is a Hopf subalgebra which is isomorphic to the $4$-dimensional Sweedler Hopf algebra. A basis of $\C$ as vector space is given by $\{a^i,ba^i,i\in\I_{0,3}\}$.

\item Denote the basis of $\mathcal{K}^{\ast}$ dual to the basis of $\mathcal{K}$ by $\{(a^i)^{\ast},(ba^i)^{\ast},i\in\I_{0,3}\}$.
From the multiplication table induced by the relations of $\mathcal{K}$, we have
\begin{align*}
\Delta(\widetilde{x})=\widetilde{x}\otimes \epsilon+\widetilde{g}\otimes \widetilde{x},\quad
\Delta(\widetilde{g})=\widetilde{g}\otimes \widetilde{g},
\end{align*}
where $\widetilde{x}=\sum_{i=0}^{3}(ba^i)\As$, $\widetilde{g}=\sum_{i=0}^{3}\xi^{-i}(a^i)\As$.
\item Let $\alpha\in G(\mathcal{K}\As)=Alg(\mathcal{K},\K)$. Since $a^{4}=1,b^2=0,ba=\xi ab$, it follows that $\alpha(a)$ is a $4$-th root of unity and $\alpha(b)=0$. Thus
    $
    G(\mathcal{K}\As)=\{\alpha_i=\sum_{j=0}^{3}\xi^{-ij}(a^j)\As\}.
  $
   Note that $\alpha_0=\epsilon$, $\alpha_i=(\alpha_1)^i$, and $G(\mathcal{K}\As)\simeq \Z_{4}$ with generator $\alpha_1$.
\item Let $\{ g^i, g^ix\}_{0\leq i<4}$ be a linear basis of $\A$. The Hopf algebra isomorphism $\phi:\A\mapsto \mathcal{K}\As$ is given by
\begin{align*}
\phi(g^i)&=\alpha_i=\sum_{j=0}^{3}\xi^{-ij}(a^j)\As,\quad
\phi(g^ix)=\theta\sum_{j=0}^{3}\xi^{-i(j+1)}(ba^j)\As,\;\text{where}\;\theta^2=2\xi.
\end{align*}
\end{enumerate}
\end{rmk}

We end up this subsection by describing explicitly the Hopf algebra structure of $\D(\mathcal{K}^{\text{cop}})$. Recall that
the Drinfeld double $\D:=\D(\mathcal{K})$ is a Hopf algebra with the tensor product coalgebra structure and the algebra structure given by
$(p\otimes a)(q\otimes b)=p\langle q_{(3)}, a_{(1)}\rangle q_{(2)}\otimes a_{(2)}\langle q_{(1)}, S^{-1}(a_{(3)})\rangle$ for $p,q\in\mathcal{K}\As,a,b\in\mathcal{K}$.
\begin{pro}
$\D:=\D(\mathcal{K}^{cop})$ as a coalgebra is isomorphic to the tensor coalgebra $\A^{cop\,op}\otimes \mathcal{K}^{cop}$, and as an algebra is generated by the elements $a, b, g, x$ satisfying the relations in $\mathcal{K}^{cop}$, the relations in $\A^{bop}$ and
\begin{align*}
ag&=ga,\quad ax-\xi xa=-\theta(ba^2-gb),\quad
bg=-gb,\quad bx-\xi xb=\theta\xi^{3}(a^{3}-ga).
\end{align*}
\end{pro}
%\begin{proof}
%Note that we have the following coproducts in $\D$
%\begin{align*}
%\De^{2}(g)&=g\otimes g\otimes g,\quad \De_{\A_2}^{2}(x)=1\otimes 1\otimes 1\otimes x+1\otimes x\otimes g+x\otimes g\otimes g,\\
%\De^2(a)&=a\otimes a\otimes a+\Lam^{-1}ba^p\otimes b\otimes a+\Lam^{-1}ba^p\otimes a^{p+1}\otimes b+\Lam^{-1}a\otimes ba^p\otimes b,\\
%\De^2(b)&=b\otimes a\otimes a+a^{p+1}\otimes b\otimes a+a^{p+1}\otimes a^{p+1}\otimes b+\Lam^{-1}b\otimes ba^p\otimes b.
%\end{align*}
%Then
%\begin{align*}
%ag&=\langle g,a\rangle ga\langle g,S(a)\rangle=ga,\\
%bg&=\langle g,a^{p+1}\rangle gb\langle g,S(a)\rangle =-gb,\\
%ax&=\Lam^{-1}\langle 1,a\rangle ba^p\langle x,S(b)\rangle+\langle 1,a\rangle xa\langle g,S(a)\rangle+\Lam^{-1}\langle x,ba^p\rangle gb\langle g,S(a)\rangle\\
%  &=\Lam^{-1}\theta\xi^{p+1}ba^p+\xi xa+\Lam^{-1}\theta\xi gb,\\
%bx&=\langle 1,a^{p+1}\rangle a^{p+1}\langle x,S(b)\rangle+\langle 1,a^{p+1}\rangle xb\langle g,S(a)\rangle+\langle x,b\rangle ga\langle g,S(a)\rangle\\
%  &=\theta\xi^{p+1}a^{p+1}+\xi xb+\theta\xi ga.
%\end{align*}
%\end{proof}

\subsection{The representation of the Drinfeld double $\D$}\label{subsecPresentation}
We describe the simple $\D$-modules. These results have been determined in \cite{GG16} and we rewrite them without the proofs since the description of $\mathcal{K}$ is different from that in \cite{GG16}.

\begin{defi}$\label{onesimple}$
Let $i\in\I_{0,3}$ and $\chi$ be an irreducible character of the cyclic group $\Z_{4}$. Denote by $\K_{\chi^{i}}$ the one-dimensional
left $\D$-module defined by
\begin{align*}
\chi^i(a)=\xi^i, \quad \chi^i(b)=0,\quad \chi^i(g)=(-1)^i,\quad \chi^i(x)=0.
\end{align*}
\end{defi}
%\begin{proof}
%Let $\chi\in G(\D\As)=\hom(\D,\K)$. Since $a^{2p}=1=g^{2p}$, we have $\chi(a)$ and $\chi(g)$ are both $2p$-th roots of unity.
%From $b^2=0$, $gb=-bg$ and $gx=-xg$, we have that $\chi(b)=\chi(x)=0$, and whence $\chi(g)^2=1$ since $x^2=1-g^2$. From
%the relation $bx-\xi xb=\theta\xi^{p+1}(a^{p+1}-ga)$, we have $\chi(a)^p=\chi(g)$. Thus $\chi$ is completely determined by $\chi(a)$. Let $\chi(a)=\xi^i$ for some $i\in Z_{2p}$, for different $i$, it is clear that these modules are pairwise non-isomorphic and any one-dimensional $\D$-module is isomorphic to $\K_{\chi^i}$ for some $i\in\Z_{2p}$.
%\end{proof}
\begin{defi}$\label{twosimple}$
For any  $(i,j)\in\Lambda=\{(i,j)\in \I_{0,3}\times \I_{0,3}\mid 2i\neq j\}$, let $V_{i,j}$ be
the $2$-dimensional left $\D$-module whose matrices defining $\D$-action with respect to a fixed basis are given by
\begin{align*}
    [a]_{i,j}&=\left(\begin{array}{ccc}
                                   \xi^i & 0\\
                                   0 &    \xi^{i+1}
                                 \end{array}\right),\quad
    [b]_{i,j}=\left(\begin{array}{ccc}
                                   0 & 1\\
                                   0 & 0
                                 \end{array}\right),\quad
    [g]_{i,j}=\left(\begin{array}{ccc}
                                   \xi^j & 0\\
                                   0 & -\xi^j
                                 \end{array}\right),\\
    [x]_{i,j}&=\left(\begin{array}{ccc}
                                   0 & \theta^{-1}\xi^{1-i}((-1)^i+\xi^j)\\
                                   \theta\xi^{i-1}((-1)^i-\xi^j) & 0
                                 \end{array}\right).
\end{align*}

\end{defi}

\begin{rmk}\label{rmkDmoddual}
For a left $\D$-module $V$, there exists a left dual module denoted by $V\As$ with the module structure given by
$(h\rightharpoonup f)(v)=f(S(h)\cdot v)$ for all $h\in \D, v\in V, f\in V\As$. A direct computation shows that $V_{i,j}\As\cong V_{-i-1,-j-2}$ for all $(i,j)\in \Lambda$.
\end{rmk}

\begin{thm}\label{thmsimplemoduleD}
There are $16$ simple $\D$-modules up to isomorphism, among which $4$ one-dimensional modules are given in Definition
$\ref{onesimple}$ and $12$ two-dimensional simple modules are given in Definition $\ref{twosimple}$.
\end{thm}
\begin{proof}
Similar to the proof of \cite[Theorem\,2.9]{GG16}.
\end{proof}

\section{Nichols algebras  in $\KYD$}\label{secNicholsalgebra}
In this section, we determine all finite-dimensional Nichols algebras in $\KYD$. As it turns out in \cite[Theorem\,4.5]{GG16} that the Nichols algebras over non-simple
indecomposable objects in $\KYD$ are infinite-dimensional, it suffices to determine all finite-dimensional Nichols algebras over the semisimple objects in $\KYD$. We first discard the Nichols algebras with infinite dimension. Then we prove that the remaining are finite-dimensional and present them by generators and relations.  It should be figured out that the Nichols algebras over the simple objects have been determined in \cite{GG16}.
\subsection{The braidings of the simple objects in $\KYD$}
We describe the braidings of the simple  objects in $\KYD$ by using the equivalence  $\KYD\cong {}_{\D}\mathcal{M}$.

\begin{pro}\label{proYDone}
Let $\K_{\chi^i}=\K\{v\}$ be a one-dimensional $\D$-module for $ i\in \I_{0,3}$. Then $\K_{\chi^i}\in\KYD$ with the Yetter-Drinfeld module structure given by
\begin{align*}
a\cdot v=\xi^i v,\quad b\cdot v=0,\quad \delta(v)=a^{2i}\otimes v.
\end{align*}
\end{pro}
\begin{proof}
Similar to the proof of \cite[Proposition\,3.1]{GG16}.
\end{proof}
%\begin{pro}
%Let $\K_{\chi^i}=\K v$ be a one-dimensional $\D$-module with $ i\in Z_6$. Then $\K_{\chi^i}\in\CYD$ with the module structure and comodule structure given by
%\begin{align*}
%a\cdot v=\xi^i v,\quad b\cdot v=0,\quad \delta(v)=a^{3i}\otimes v.
%\end{align*}
%\end{pro}
%\begin{proof}
%Since $\K_{\chi^i}=\K\{v\}$ is a one dimensional $\D$-module with $ i\in\I_{0,3}$, the $\C$-action must be given by the restriction of the character of $\D$ given by lemma $\ref{onesimple}$ and the coaction must be of the form $\delta(v)=h\otimes v$ where $h\in G(\mathcal{K})=\{1,a^2\}$ such that $\langle g, h\rangle v=(-1)^iv$. It follows that the action is given by $a\cdot v=\xi^i v,\, b\cdot v=0$ and the coaction is given by $\delta(v)=a^{2i}\otimes v$.
%\end{proof}
\begin{rmk}\label{rmkVA1}
Let $\K_{\chi^i}=\K\{v\}\in\KYD$ for $ i\in \I_{0,3}$. Then $\K_{\chi^i}\in{}_{\A}^{\A}\mathcal{YD}$ with the Yetter-Drinfeld module structure given by
\begin{align*}
g\cdot v=(-1)^iv,\quad x\cdot v=0,\quad \delta(v)=g^i\otimes v.
\end{align*}
\end{rmk}
\begin{pro}
Let $V_{i,j}=\K\{v_1,v_2\}$ be a two-dimensional simple $\D$-module for $(i,j)\in\Lambda$. Then $V_{i,j}\in\KYD$ with the Yetter-Drinfeld module structure  given by
\begin{gather*}
a\cdot v_1=\xi^iv_1,\quad b\cdot v_1=0,\quad a\cdot v_2=\xi^{i+1} v_2,\quad b\cdot v_2=v_1,\\
%\end{align*}
%and the comodule structure given by
%\begin{align*}
\delta(v_1)=a^{-j}\otimes v_1+\xi^{i-1}((-1)^i-\xi^j)ba^{-1-j}\otimes v_2,\;
\delta(v_2)=a^{2-j}\otimes v_2+\frac{1}{2}\xi^{-i}((-1)^i+\xi^j)ba^{1-j}\otimes v_1.
\end{gather*}
\end{pro}
\begin{proof}
Similar to the proof of \cite[Proposition\,3.3]{GG16}.
\end{proof}

\begin{rmk}\label{rmkVA2}
Let $V_{i,j}=\K\{v_1,v_2\}\in\KYD$ for $(i,j)\in\Lambda$. Then $V_{i,j}\in\AYD$ with the Yetter-Drinfeld module structure given by
\begin{align*}
g\cdot v_1=\xi^{-j}v_1,\quad x\cdot v_1=x_2\xi^{2-j}v_2,\quad g\cdot v_2=\xi^{2-j}v_2,\quad x\cdot v_2=x_1\xi^{-j}v_1,\\
\delta(v_1)=g^i\otimes v_1,\quad \delta(v_2)=g^{i+1}\otimes v_2+\theta^{-1}(-1)^{i+1}\xi^{-i-1}g^ix\otimes v_1,
\end{align*}
where $x_1=\theta^{-1}\xi^{1-i}((-1)^i+\xi^j)$ and $x_2=\theta\xi^{i-1}((-1)^i-\xi^j)$.
\end{rmk}
%\begin{pro}
%Let $\Pp_j$ be the projective cover of the one-dimensional $\D$-module $\K_{\chi^j}$ for $j\in\I_{0,3}$. Then $\Pp_j\in\KYD $ with its module structure given by $\eqref{eqprojective0}$ and $\eqref{eqprojectivei}$ and its comodule structure given by
%\begin{align*}
%\delta(p_{1,j})&=(a^2)^j\otimes p_{1,j}+\theta^{-1}(-1)^jba^{-1}(a^2)^j\otimes (\theta p_{2,j}+2p_{3,j}),\\
%\delta(p_{2,j})&=(a^2)^{j+1}\otimes p_{2,j}+(-1)^jba^{-1}a^{2(j+1)}\otimes p_{4,j},\\
%\delta(p_{3,j})&=(a^2)^{j+1}\otimes p_{3,j}+\xi^{-1}\theta^{-1}(-1)^jba^{-1}a^{2(j+1)}\otimes p_{4,j},\\
%\delta(p_{4,j})&=(a^2)^j\otimes p_{4,j}.
%\end{align*}
%\end{pro}
%
%
Using the braiding formula $\eqref{equbraidingYDcat}$ in $\KYD$,  we describe the braidings of the simple objects in $\KYD$.
\begin{pro}\label{braidingone}
Let $\K_{\chi^i}=\K\{v\}\in\KYD$ for $i\in\I_{0,3}$. Then the braiding of $\K_{\chi^i}$ is $c(v\otimes v)=(-1)^iv\otimes v$.
\end{pro}
\begin{pro}\label{probraidsimpletwo}
Let $V_{i,j}=\K\{v_1,v_2\}\in\KYD$ for $(i,j)\in\Lambda$. Then the braiding of $V_{i,j}$ is given by
%\begin{enumerate}
%  \item If $j=0$, $i\in\{1,3,5\}$,
  \begin{align*}
   c(\left[\begin{array}{ccc} v_1\\v_2\end{array}\right]\otimes\left[\begin{array}{ccc} v_1~v_2\end{array}\right])=
   \left[\begin{array}{ccc}
   \xi^{-ij}v_1\otimes v_1    & \xi^{-j(i+1)}v_2\otimes v_1+[\xi^{-ij}+\xi^{(i+1)(2-j)}]v_1\otimes v_2\\
   (-1)^i\xi^{-ij}v_1\otimes v_2   & \xi^{(i+1)(2-j)}v_2\otimes v_2+\frac{1}{2}\xi^{1-ij-j}[(-1)^i+\xi^j]v_1\otimes v_1
         \end{array}\right].
  \end{align*}
\end{pro}
%
%
%\begin{pro}\label{probraidingprocover}
%The braiding of $\Pp_j$ is given by
%\begin{align*}
%   c(p_{1,j}\otimes\left[\begin{array}{ccc} p_{1,j}\\p_{2,j}\\p_{3,j}\\p_{4,j}\end{array}\right])&=
%   \left[\begin{array}{ccc} (-1)^jp_{1,j} \\p_{2,j}\\p_{3,j}\\(-1)^jp_{4,j}\end{array}\right]\otimes p_{1,j}+
%   \left[\begin{array}{ccc}p_{3,j} \\ (-1)^j\theta^{-1}\xi^{-1}p_{4,j}\\0\\0\end{array}\right]\otimes (\theta p_{2,j}{+}2p_{3,j}),\\
%   c(p_{2,j}\otimes\left[\begin{array}{ccc} p_{1,j}\\p_{2,j}\\p_{3,j}\\p_{4,j}\end{array}\right])&=
%   \left[\begin{array}{ccc}p_{1,j} \\(-1)^{j+1}p_{2,j}\\(-1)^{j+1}p_{3,j}\\p_{4,j}\end{array}\right]\otimes p_{2,j}+
%   \left[\begin{array}{ccc}2(-1)^{j+1}\xi^{-1}\theta^{-1} p_{3,j}\\2\xi^{p-2}\theta^{-2}(-1)^jp_{4,j}\\0\\0\end{array}\right]\otimes p_{4,j},\\
%   c(p_{3,j}\otimes\left[\begin{array}{ccc} p_{1,j}\\p_{2,j}\\p_{3,j}\\p_{4,j}\end{array}\right])&=
%   \left[\begin{array}{ccc}p_{1,j} \\(-1)^{j+1}p_{2,j}\\(-1)^{j+1}p_{3,j}\\p_{4,j}\end{array}\right]\otimes p_{3,j}+
%   \left[\begin{array}{ccc}(-1)^{j}\xi^{-1} p_{3,j}\\\theta^{-1}\xi^{p-2} p_{4,j}\\0\\0\end{array}\right]\otimes p_{4,j},\\
%   c(p_{4,j}\otimes\left[\begin{array}{ccc} p_{1,j}\\p_{2,j}\\p_{3,j}\\p_{4,j}\end{array}\right])&=
%   \left[\begin{array}{ccc}(-1)^jp_{1,j} \\p_{2,j}\\p_{3,j}\\(-1)^jp_{4,j}\end{array}\right]\otimes p_{4,j}.
%\end{align*}
%\end{pro}
%
%
%

\subsection{Nichols algebras over the simple objects in $\KYD$}
We study the Nichols algebras over the simple objects in $\KYD$.
\subsubsection{Nichols algebras over the one-dimensional objects in $\KYD$}
\begin{lem}\label{lemNicholbyone}
The Nichols algebra $\BN(\K_{\chi^k})$ over $\K_{\chi^k}=\K\{v\}$ for $k\in\I_{0,3}$ is
\begin{align*}
\BN(\K_{\chi^k})=\begin{cases}
\K[v] &\text{~if~} k\in\{0,2\};\\
\bigwedge \K_{\chi^k} &\text{~if~} k\in\{1,3\}.\\
\end{cases}
\end{align*}
Moreover, let $V=\oplus_{i\in I}V_i$, where $V_i\cong \K_{\chi^{k_i}}$ with $k_i\in\{1,3\}$, and $I$ is a finite index set.
Then $\BN(V)=\bigwedge V\cong \otimes_{i\in I}\BN(V_i)$.
\end{lem}
\begin{proof}
By Propositions \ref{proYDone} \& \ref{braidingone}, the braiding $c=-\tau$, where $\tau$ is a flip. Thus the lemma follows.
\end{proof}

 %Recall that $\gr\A=\BN(W)\sharp \K[\Gamma]$, where $\Gamma\cong \Z_4$ with the generator $g$ and $W:=\K\{x\}\in{}_{\Gamma}^{\Gamma}\mathcal{YD}$ with the Yetter-Drinfeld module structure given by $g\cdot x=-x$ and $\delta(x)=g\otimes x$.
\subsubsection{Nichols algebras over the two-dimensional simple objects in $\KYD$}
Let $V_{i,j}=\K\{v_1,v_2\}\in\KYD$ with $(i,j)\in\Lambda$. Then $V_{i,j}\in\AYD$ with the Yetter-Drinfeld module structure given  by Remark \ref{rmkVA2}. Let $B_{i,j}$  be the subalgebra of $\BN(V_{i,j})\sharp\A$ generated by $g,x,v_1$. Then $B_{i,j}$ is a pointed Hopf algebra with $\G(B_{i,j})\cong\Gamma$, which is isomorphic to the quotient of  $\BN(X_{i,j})\sharp\K[\Gamma]$ by the relation $x^2=1-g^2$, where $X_{i,j}=\K\{x,v_1\}\in{}_{\Gamma}^{\Gamma}\mathcal{YD}$ with the Yetter-Drinfeld module structure given by
\begin{align*}
g\cdot x&=-x,\quad  g\cdot v_1=\xi^{-j}v_1,\quad \delta(x)=g\otimes x,\quad\delta(v_1)=g^i\otimes v_1.
\end{align*}
It is easy to see that $\gr B_{i,j}\cong \BN(X_{i,j})\sharp\K[\Gamma]$ and $\BN(X_{i,j})$ is of diagonal type with the generalized Dynkin diagram  given by \Dchaintwo{-1}{\xi^{2i-j}}{\xi^{-ij}}.
\begin{lem}\label{lemOne}
Let $\Lambda\As=\{(i,j)\in\Lambda\mid ij=0\text{~or~}(i+1)(2-j)=0 \text{~mod~} 4\}$. Then $\dim\,\BN(V_{k,l})=\infty$ for all $(k,l)\in\Lambda\As$.
\end{lem}
\begin{proof}
Let $(i,j)\in\Lambda\As$. If $ij=0\text{~mod~}4$, i.~e.~$\xi^{-ij}=1$, then $c(v_1\otimes v_1)=v_1\otimes v_1$ and whence $\dim\BN(X_{i,j})=\infty$. If $ij\neq 0\text{~mod~}4$, then $(i,j)\in\{(3,1),(3,3)\}$. We apply the reflection described as follows:
$$
\Dchaintwo{-1}{\xi^{\pm 1}}{\xi^{\pm 1}} \xymatrix{\ar@/^0.15pc/@{->}[r]^{x} & }
\xymatrix@C-4pt{\overset{-1}{\underset{\  }{\circ}}\ar
@{-}[r]^{\xi^{\mp 1}}
& \overset{1}{\underset{\  }{\circ}}}
$$
It follows that $\dim\BN(X_{i,j})=\infty$. As $\dim\BN(V_{i,j})\sharp\A\geq \dim B_{i,j}=\dim\text{gr}\,B_{i,j}=\BN(X_{i,j})\sharp\K[\Gamma]$, it follows that $\dim \BN(V_{i,j})=\infty$ for $(i,j)\in\Lambda\As$.
\end{proof}

Note that $\Lambda-\Lambda\As=\{(1,1),(1,3),(2,1),(2,3)\}$. We shall show that $\dim\BN(V_{i,j})<\infty$ for $(i,j)\in\Lambda-\Lambda\As$.
\begin{pro}\label{proNicholsTwoOfK}
Let $V$ be a simple object in $\KYD$ such that the Nichols algebra $\BN(V)$ is finite-dimensional. Then $V$ is isomorphic either to $\K_{\chi}$, $\K_{\chi^{-1}}$, $V_{1,1}$, $V_{1,3}$, $V_{2,1}$ or $V_{2,3}$. Moreover, the generators and relations of $\BN(V)$ are given by
\center
%\caption{Finite-dimensional Nichols algebras}
\begin{tabular}{|c|c|c|}
  \hline
  $V$           &   relations of $\BN(V)$ with generators $v_1,v_2$         & $\dim\BN(V)$ \\\hline
$\K_{\chi^1},l=1,3$   &   $v_1^2=0$                                       &   2   \\\hline
$V_{1,j},j=1,3$       &   $v_1^2+2\xi v_2^2=0, v_1v_2+v_2v_1=0, v_1^4=0  $    &  8  \\\hline
$V_{2,j},j=1,3$       &    $ v_1^2=0, v_2v_1+\xi^jv_1v_2=0, v_2^4=0$   &  8  \\\hline
\end{tabular}
\end{pro}
\begin{proof}
Using the braidings in Proposition \ref{probraidsimpletwo}, the lemma follows by \cite[Proposition\; 3.10 \& 3.11]{AGi17}.
\end{proof}
\begin{rmk}\label{rmksplitNichols}
Let $(i,j)\in\Lambda-\Lambda\As$. We claim that $ B_{i,j}\cong \BN(V_{i,j})\sharp\A$ as Hopf algebras.
Recall that $B_{i,j}$ is a Hopf subalgebra of $\BN(V_{i,j})\sharp\A$. It suffices to show that $\dim B_{i,j}=\dim \BN(V_{i,j})\sharp\A$. Indeed, if $i=1$ or $2$, then  the Dynkin diagram of $\BN(X_{i,j})$ is \Dchaintwo{-1}{\xi^j}{\xi^{-j}} or \Dchaintwo{-1}{\xi^{-j}}{-1}, respectively.
Then $\dim\BN(X_{i,j})=8$ by \cite{An13,An15} and whence $\dim B_{i,j}=\dim\gr B_{i,j}=64=\dim\BN(V_{i,j})\sharp\A$.

We claim that $\BN(V_{i,j})\sharp\gr \A\cong\BN(X_{i,j})\sharp\K[\Gamma]$.
Recall that $\A$ is a cocycle deformation of $\gr A\cong\BN(W)\sharp\K[\Gamma]$, that is, $\A\cong(\gr\A)_{\sigma}$ for some Hopf 2-cocycle $\sigma$. By \cite[Example\,5.1]{GM10},
$$\sigma=\epsilon\otimes\epsilon+\zeta, \text{ where } \zeta(x^mg^n,x^kg^l)=(-1)^{nk}\delta_{2,m+k}.$$
By \cite[Theorem\,2.7]{MO99}, a direct computation shows that $V_{i,j}=\K\{v_1,v_2\}\in\grAYD$ by
\begin{gather*}
g\cdot v_1=\xi^{-j}v_1,\quad x\cdot v_1=A_1v_2+A_2v_1,\quad g\cdot v_2=\xi^{2-j}v_2,\quad x\cdot v_2=B_1v_1+B_2v_2,\\
\delta(v_1)=g^i\otimes v_1,\quad \delta(v_2)=g^{i+1}\otimes v_2+\theta^{-1}(-1)^{i+1}\xi^{-i-1}g^ix\otimes v_1,
\end{gather*}
where $A_1,A_2,B_1,B_2\in\K$.
Then by \cite[Proposition\,8.8]{HS13}, $\BN(V_{i,j})\sharp \BN(W)\cong \BN(X_{i,j})$. Consequently, the claim follows.

Note that $V_{i,j}\in\KYD$ $($or $\grAYD$$)$ is characterized by $(i,j)\in\Lambda$ and $V_{i,j}\cong V_{p,q}$ in $\grAYD$ if and only if $V_{i,j}\cong V_{p,q}$ in $\KYD$ for $(i,j),(p,q)\in\Lambda$.  From \cite[Proposition\,8.6]{HS13} and the preceding discussion, the Nichols algebras of dimension bigger than $2$ in Proposition \ref{proNicholsTwoOfK} can be obtained $($up to isomorphism$)$ by splitting the Nichols algebras of diagonal type.
\end{rmk}

\subsection{Nichols algebras over the semisimple objects in $\KYD$}
We determine all finite-dimensional Nichols algebras over semisimple objects and present them by generators and relations.
\subsubsection{Nichols algebras over $V_{i,j}\oplus\K_{\chi^k}$ in $\KYD$}
Let $V_{i,j}=\K\{v_1,v_2\}\in\KYD$ for $(i,j)\in\Lambda-\Lambda\As$ and $\K_{\chi^k}=\K\{v_3\}\in\KYD$ for $k\in\{1,3\}$. Then $V_{i,j}\oplus\K_{\chi^k}\in\KYD$ with the Yetter-Drinfeld module structure given by
\begin{align*}
a\cdot v_1&=\xi^i v_1,\quad b\cdot v_1=0,\quad \delta(v_1)=a^{-j}\otimes v_1+\xi^{i-1}((-1)^i-\xi^j)ba^{-1-j}\otimes v_2,\\
a\cdot v_2&=\xi^{i+1} v_2,\quad b\cdot v_2=v_1,\quad \delta(v_2)=a^{2-j}\otimes v_2+\frac{1}{2}\xi^{-i}((-1)^i+\xi^j)ba^{1-j}\otimes v_1,\\
a\cdot v_3&=\xi^k v_3,\quad b\cdot v_3=0,\quad \delta(v_3)=a^{2}\otimes v_3.
\end{align*}
Moreover, the braiding of $V_{i,j}\oplus\K_{\chi^k}$ is given as follows: $c(\left[\begin{array}{ccc} v_1\\v_2\\v_3\end{array}\right]\otimes\left[\begin{array}{ccc} v_1~v_2~v_3\end{array}\right])=$
\begin{align}\label{BraidingTwoone}
\left[\begin{array}{ccc}
   \xi^{-ij}v_1\otimes v_1    & \xi^{-j(i+1)}v_2\otimes v_1+[\xi^{-ij}+\xi^{(i+1)(2-j)}]v_1\otimes v_2 & \xi^{-kj}v_3\otimes v_1\\
   (-1)^i\xi^{-ij}v_1\otimes v_2   & \xi^{(i+1)(2-j)}v_2\otimes v_2+\frac{1}{2}\xi^{1-ij-j}[(-1)^i+\xi^j]v_1\otimes v_1&\xi^{(2-j)k}v_3\otimes v_2\\
   (-1)^{ik}v_1\otimes v_3 & (-1)^{(i+1)k}v_2\otimes v_3 & (-1)^k v_3\otimes v_3
         \end{array}\right].
\end{align}
 In particular, the braided vector subspace $\K\{v_1,v_3\}$ is of diagonal type.

Note that $V_{i,j}\oplus\K_{\chi^k}\in\AYD$  with the Yetter-Drinfeld module structure given by Remarks \ref{rmkVA1} \& \ref{rmkVA2}. Let $C_{i,j,k}$ be the subalgebra of $\BN(V_{i,j}\oplus\K_{\chi^k})\sharp\A$ generated by $v_1,v_3,g,x$.  Then  $C_{i,j,k}$ is a pointed Hopf algebra with $\G(C_{i,j,k})\cong\Gamma$, which is isomorphic to the quotient of  $\BN(Y_{i,j,k})\sharp\K[\Gamma]$ by the relation $x^2=1-g^2$, where $Y_{i,j,k}=\K\{x,v_1,v_3\}\in{}_{\Gamma}^{\Gamma}\mathcal{YD}$ with the Yetter-Drinfeld module structure given by
\begin{align*}
g\cdot x&=-x,\quad  g\cdot v_1=\xi^{-j}v_1,\quad g\cdot v_3=-v_3\quad \delta(x)=g\otimes x,\quad\delta(v_1)=g^i\otimes v_1,\quad\delta(v_3)=g^k\otimes v_3.
\end{align*}
It is easy to see that $\text{gr}\,C_{i,j,k}\cong \BN(Y_{i,j,k})\sharp\K[\Gamma]$ and $\BN(Y_{i,j,k})$ is of diagonal type with the generalized Dynkin diagram given by \xymatrix@C+9pt{
\overset{-1}{\underset{x}{\circ}}\ar  @ {-}[r]^{(-1)^i\xi^{-j}}  & \overset{\xi^{-ij}}{\underset{v_1
}{\circ}}\ar  @{-}[r]^{(-)^i\xi^{-kj}}
& \overset{-1}{\underset{v_3}{\circ}} }.
\begin{lem}\label{lemTwoOne}
Let $j\in\{1,3\}$ and $V,\;W$ be two simple objects in $\KYD$. If $V\oplus W$ is isomorphic either to $V_{1,j}\oplus\K_{\chi^3}$ or $V_{2,j}\oplus\K_{\chi}$, then $\BN(V\oplus W)$ is infinite-dimensional.
\end{lem}
\begin{proof}
%\begin{enumerate}
The Dynkin diagram of $Y_{1,j,3}$ or $Y_{2,j,1}$ is \Dchainthree{-1}{x}{\xi^j}{\xi^{-j}}{v_1}{\xi^{-j}}{-1}{v_3} or \Dchainthree{-1}{x}{\xi^{-j}}{-1}{v_1}{\xi^{-j}}{-1}{v_3}, respectively. We apply the reflection at several points  described as follows,
  $$
   \Dchainthree{-1}{x}{\xi^{-j}}{-1}{v_1}{\xi^{-j}}{-1}{v_3} \xymatrix{\ar@/^0.15pc/@{->}[r]^{x} & }
   \Dchainthree{-1}{1 }{\xi^j}{\xi^{-j}}{2 }{\xi^{-j}}{-1}{3 } \xymatrix{\ar@/^0.15pc/@{->}[r]^{3} & }
   \Dchainthree{-1}{1 }{\xi^j}{1}{2 }{\xi^{j}}{-1}{3 }.
  $$
By Theorem \ref{thmReflection}, the lemma is proved.
\end{proof}

\begin{pro}\label{proNicholsTwoone}
Let $V_{i,j}$ be a two-dimensional simple object and $\K_{\chi^k}$ be a one-dimensional simple object in $\KYD$ such that $\dim\BN(V_{i,j}\oplus\K_{\chi^k})<\infty$ for $(i,j)\in\Lambda,k\in\I_{0,3}$. Then $V_{i,j}\oplus\K_{\chi^k}$ is isomorphic either to $V_{1,1}\oplus \K_{\chi}$, $V_{1,3}\oplus\K_{\chi}$, $V_{2,1}\oplus\K_{\chi^3}$ or $V_{2,3}\oplus\K_{\chi^{3}}$. Moreover, $\BN(V_{1,j}\oplus\K_{\chi})$ for $j=1,3$ is generated by the elements $v_1$, $v_2$, $v_3$ such that $v_1,v_2$ satisfy the relations of $\BN(V_{1,j})$, $v_3^2=0$ and all together they satisfy the relations:
\begin{gather}
 v_3v_1^2+(1-\xi^j)v_1v_3v_1-\xi^j v_1^2v_3=0,\quad\xi^j v_1v_2v_3-\xi^j v_1v_3v_2+v_2v_3v_1+v_3v_1v_2=0,\label{eqV1j1-1}\\
\frac{1}{2}\xi(1+\xi^{-j})(v_1v_3)^2(v_2v_3)^2+(v_2v_3)^4+(v_3v_2)^4=0.\label{eqV1j1-2}
\end{gather}
$\BN(V_{2,j}\oplus\K_{\chi^3})$ for $j=1,3$ is generated by the elements $v_1$, $v_2$, $v_3$ such that $v_1,v_2$ satisfy the relations of $\BN(V_{2,j})$, $v_3^2=0$ and all together they satisfy the relations:
\begin{gather}
 v_1v_2v_3+v_1v_3v_2+v_2v_3v_1+v_3v_2v_1=0,\quad (v_3v_1)^4+(v_1v_3)^4=0,\\
  \frac{1}{2}\xi(1+\xi^{-j})v_1v_3v_1-\xi^j v_2^2v_3+(1-\xi^j)v_2v_3v_2+v_3v_2^2=0.
\end{gather}
\end{pro}
\begin{proof}
We prove the statement only for $\BN(V_{1,j}\oplus \K_{\chi})$, being the proof for $\BN(V_{2,j}\oplus \K_{\chi^3})$ completely analogous.
Using the braiding \eqref{BraidingTwoone}, a tedious computation shows that relations \eqref{eqV1j1-1} \eqref{eqV1j1-2}  are annihilated by $\partial_i$, $i=1,2,3$.  Hence the quotient $\BN$ of $\BN(V_{1,j})\otimes\BN(\K_{\chi})$  by \eqref{eqV1j1-1} \eqref{eqV1j1-2}  projects onto $\BN(V_{1,j}\oplus \K_{\chi}))$.

Let $B:=\{v_1^{n_1}v_2^{n_2}(v_3v_1)^{n_{31}}(v_3v_2)^{n_{32}}v_3^{n_3}:\;n_1,n_{32}\in\I_{0,3},n_2,n_{31},n_3\in\I_{0,1}\}$. Note that $|B|=128$. We claim that  the subspace $I$ linearly spanned by $B$ is a left ideal of $\BN$.  Then $B$ linearly  generates $\mathfrak{B}$ since $1\in I$. Indeed, for this, it is enough to show that $gI\subset I$ for $g\in\{v_1,v_2,v_3\}$ and they can be obtained easily from the defining relations of $\BN$ and additional relations  induced by the defining relations:
\begin{align*}
(v_3v_1)^2+\xi^j(v_1v_3)^2=0,\quad v_3v_2v_3v_1-\xi^jv_3v_1v_3v_2-\xi^jv_2v_3v_1v_3-v_1v_3v_2v_3.
\end{align*}

We claim that $\dim\BN(V_{1,j}\oplus \K_{\chi})\geq 128=|B|$. Indeed, a direct computation shows that the Dynkin diagram of $\BN(Y_{1,j,1})$ is \Dchainthree{-1}{x}{\xi^j}{\xi^{-j}}{v_1}{\xi^{j}}{-1}{v_3} and whence $\dim\BN(Y_{1,j,1})=256$ by \cite{An13,An15}. It follows that $\dim\BN(V_{1,j}\oplus \K_{\chi})\sharp\mathcal{K}=\dim\BN(V_{1,j}\oplus \K_{\chi})\sharp\A \geq\dim \BN(Y_{1,j,1})\sharp\K[\Gamma]=1024$. Since $\dim\mathcal{K}=8$, the claim follows. Consequently, $\BN\cong\BN(V_{1,j}\oplus \K_{\chi})$.
\end{proof}

\begin{rmk}\label{rmksplitNichols1}
Let $(i,j,k)\in\{(1,1,1),(1,3,1),(2,1,3),(2,3,3)\}$. Similar to Remark \ref{rmksplitNichols}, $ C_{i,j,k}\cong \BN(V_{i,j}\oplus\K_{\chi^k})\sharp\A$ and $\BN(V_{i,j}\oplus\K_{\chi^k})\sharp\gr\A\cong\BN(Y_{i,j,k})\sharp\K[\Gamma]$. If $i=1$ or $2$, then the Dynkin diagram of $\BN(Y_{i,j,k})$ is \Dchainthree{-1}{x}{\xi^j}{\xi^{-j}}{v_1}{\xi^{j}}{-1}{v_3} or \Dchainthree{-1}{x}{\xi^{-j}}{-1}{v_1}{\xi^{j}}{-1}{v_3}, respectively. Moreover, $\BN(V_{i,j}\oplus\K_{\chi^k})\sharp \BN(W)\cong  \BN(Y_{i,j,k})$ and the Nichols algebras in Proposition \ref{proNicholsTwoone} can be obtained by splitting the Nichols algebras of diagonal type.
\end{rmk}
\subsubsection{Nichols algebras over $V_{i,j}\oplus V_{k,\ell}$ in $\KYD$}
Let $V_{i,j}=\K\{v_1,v_2\}$ and $V_{k,\ell}=\K\{e_1,e_2\}$ be two simple objects in $\KYD$ for $(i,j),(k,\ell)\in\Lambda-\Lambda\As$.
 Then $V_{i,j}\oplus V_{k,\ell}\in\KYD$ with the Yetter-Drinfeld module structure given by
\begin{align*}
a\cdot v_1&=\xi^i v_1,\quad b\cdot v_1=0,\quad \delta(v_1)=a^{-j}\otimes v_1+\xi^{i-1}((-1)^i-\xi^j)ba^{-1-j}\otimes v_2,\\
a\cdot v_2&=\xi^{i+1} v_2,\quad b\cdot v_2=v_1,\quad \delta(v_2)=a^{2-j}\otimes v_2+\frac{1}{2}\xi^{-i}((-1)^i+\xi^j)ba^{1-j}\otimes v_1,\\
a\cdot e_1&=\xi^k v_1,\quad b\cdot e_1=0,\quad \delta(e_1)=a^{-\ell}\otimes v_1+\xi^{k-1}((-1)^k-\xi^{\ell})ba^{-1-\ell}\otimes v_2,\\
a\cdot e_2&=\xi^{k+1} v_2,\quad b\cdot v_2=v_1,\quad \delta(v_2)=a^{2-\ell}\otimes v_2+\frac{1}{2}\xi^{-k}((-1)^k+\xi^{\ell})ba^{1-\ell}\otimes v_1.
\end{align*}
Note that $V_{i,j}\oplus V_{k,\ell}\in\AYD$ with the Yetter-Drinfeld module structure given by Remark \ref{rmkVA2}. Let $D_{i,j,k,\ell}$ be the subalgebra of $\BN(V_{i,j}\oplus V_{k,\ell})\sharp\A$ generated by $v_1,e_1,g,x$.  Then $D_{i,j,k,\ell}$ is a pointed Hopf algebra with $G(D_{i,j,k,\ell})\cong\Gamma$, which is isomorphic to the quotient of  $\BN(Z_{i,j,k,\ell})\sharp\K[\Gamma]$ by the relation $x^2=1-g^2$, where $Z_{i,j,k,\ell}=\K\{x,e_1,v_1\}\in{}_{\Gamma}^{\Gamma}\mathcal{YD}$ with the Yetter-Drinfeld module structure given by
\begin{align*}
g\cdot x&=-x,\quad  g\cdot v_1=\xi^{-j}v_1,\quad g\cdot e_1=\xi^{-\ell}e_1\quad \delta(x)=g\otimes x,\quad\delta(v_1)=g^i\otimes v_1,\quad\delta(e_1)=g^k\otimes e_1.
\end{align*}
It is easy to see that $\text{gr}\,D_{i,j,k,\ell}\cong \BN(Z_{i,j,k,\ell})\sharp\K[\Gamma]$ and $\BN(Z_{i,j,k,\ell})$ is of diagonal type with the generalized Dynkin diagram given by \Dtriangle{-1}{(-1)^k\xi^{-\ell}}{(-1)^i\xi^{-j}}{\xi^{-k\ell}}{\xi^{-kj-i\ell}}{\xi^{-ij}.}

\begin{lem}
Let $j\in\{1,3\}$ and $V,\;W$ be two simple objects in $\KYD$. Then $\BN(V\oplus W)$ is infinite-dimensional if $V\oplus W$ is isomorphic either to $V_{1,j}\oplus V_{1,j}$, $V_{2,j}\oplus V_{2,j}$, $V_{2,j}\oplus V_{1,j}$ or $V_{2,j}\oplus V_{1,-j}$.
\end{lem}
\begin{proof}
\begin{enumerate}
  \item For $V_{2,j}\oplus V_{2,j}$, the Dynkin diagram is  \Dchainthree{-1}{v_1}{\xi^{-j}}{-1}{x}{\xi^{-j}}{-1}{e_1}. Then as shown in the proof in Lemma \ref{lemTwoOne},  $\dim\BN(V_{2,j}\oplus V_{2,j})=\infty$.
  \item For $V_{1,j}\oplus V_{1,j}$ or $V_{2,j}\oplus V_{1,-j}$, the Dynkin diagram is
  \begin{align*}
  \xymatrix@R-12pt{  &\overset{-1}{\underset{x}{\circ}} \ar  @{-}[dl]_{\xi^j}\ar  @{-}[dr]^{\xi^{j}} & \\
\overset{\xi^{-j}}{\underset{e_1}{\circ}} \ar  @{-}[rr]^{-1}  &  &\overset{\xi^{-j}}{\underset{v_1}{\circ}},\text{~~or~~} }
\xymatrix@R-12pt{  &\overset{-1}{\underset{x}{\circ}} \ar  @{-}[dl]_{\xi^{-j}}\ar  @{-}[dr]^{\xi^{-j}} & \\
\overset{\xi^{j}}{\underset{e_1}{\circ}} \ar  @{-}[rr]^{\xi^j}  &  &\overset{-1}{\underset{v_1}{\circ}}.}
  \end{align*}
 Since \xymatrix@C-4pt{\overset{\xi^{-j}}{\underset{e_1 }{\circ}}\ar
@{-}[r]^{-1}
& \overset{\xi^{-j}}{\underset{ e_2}{\circ}}} is of affine Cartan type and \xymatrix@C-4pt{\overset{\xi^{j}}{\underset{e_1 }{\circ}}\ar
@{-}[r]^{\xi^j}
& \overset{-1}{\underset{ e_2}{\circ}}} has an infinite system as shown in the proof of Lemma \ref{lemTwoOne}, it follows that $\dim\BN(V_{1,j}\oplus V_{1,j})=\infty= \dim\BN(V_{2,j}\oplus V_{1,-j})$.

  \item For $V_{2,j}\oplus V_{1,j}$, the Dynkin diagram is  \Dtriangle{-1}{\xi^j}{\xi^{-j}}{\xi^{-j}}{\xi^j}{-1}.

Then $\dim\BN(V_{2,j}\oplus V_{1,j})=\infty$ by \cite[Lemma\,9.(ii)]{H09}.
\end{enumerate}

\end{proof}

\begin{pro}\label{proTwoTwo}
Let $V_{i,j}$ and $V_{k,\ell}$ be two simple objects in $\KYD$ for $(i,j),(k,\ell)\in\Lambda$ such that $\dim\BN(V_{i,j}\oplus V_{k,\ell})<\infty$. Then $V_{i,j}\oplus V_{k,\ell}$ is isomorphic either to  $V_{1,1}\oplus V_{1,3}$ or $V_{2,1}\oplus V_{2,3}$. Moreover, $\BN(V_{i,1}\oplus V_{i,3})$ for $i\in\{1,2\}$ is  generated by elements $v_1,v_2,e_1,e_2$ such that $v_1,v_2$ satisfy the relations in $\BN(V_{i,1})$, $e_1,e_2$ satisfy the relations in $\BN(V_{i,3})$ and all together they satisfy the relations:
\begin{align}\label{eqV1113}
e_1v_1-\xi v_1e_1=0,\  e_2v_1+e_1v_2+v_2e_1+v_1e_2=0,\  e_2v_2+v_2e_2+\frac{1}{2}(1-\xi)v_1e_1=0,\ \text{for $i=1$};\\
e_1v_1+ v_1e_1=0,\  e_2v_1-e_1v_2-\xi v_2e_1+\xi v_1e_2=0,\  e_2v_2-\xi v_2e_2+\frac{1}{2}(\xi-1)v_1e_1=0,\quad\text{for $i=2$}.
\end{align}
In particular, $\dim\BN(V_{1,1}\oplus V_{1,3})=128=\dim\BN(V_{2,1}\oplus V_{2,3})$.
\end{pro}
\begin{proof}
We prove the statement only for $\BN(V_{1,1}\oplus V_{1,3})$, being the proof for $\BN(V_{2,1}\oplus V_{2,3})$ completely analogous. A straight computation shows that relations \eqref{eqV1113} are primitive in $T(V_{1,1}\oplus V_{1,3})$. Hence the quotient $\BN$ of $\BN(V_{1,1})\otimes\BN(V_{1,3})$ by \eqref{eqV1113} projects onto $\BN(V_{1,1}\oplus V_{1,3})$.

Let $B:=\{v_1^{n_1}v_2^{n_2}(e_1v_2)^{n_{32}}e_1^{n_3}e_2^{n_4}:\;n_1,n_3\in\I_{0,3},n_2,n_{32},n_3\in\I_{0,1}\}$. Note that $|B|=128$. We claim that  the subspace $I$ linearly spanned by $B$ is a left ideal of $\BN$.  Then $B$ linearly  generates $\mathfrak{B}$ since $1\in I$. Indeed, for this, it is enough to show that $gI\subset I$ for $g\in\{v_1,v_2,e_1,e_2\}$ and they can be obtained easily from the defining relations of $\BN$ and additional relations  induced by the defining relations:
\begin{align*}
e_1^2v_2+(1+\xi)e_1v_2e_1+\xi v_2e_1^2,\quad (e_1v_2)^2-\xi (v_2e_1)^2+\frac{1}{2}(1+\xi)v_1^2e_1^2.
\end{align*}

We claim that $\dim\BN(V_{1,1}\oplus V_{1,3})\geq 128=|B|$. Indeed, a direct computation shows that the Dynkin diagram of $\BN(Z_{1,1,1,3})$ is \Dchainthree{\xi^{-1}}{v_1}{\xi}{-1}{x}{\xi^{-1}}{\xi}{e_1} and whence $\dim\BN(Z_{1,1,1,3})=256$ by \cite{An13,An15}. It follows that $\dim\BN(V_{1,1}\oplus V_{1,3})\sharp\mathcal{K}=\dim\BN(V_{1,1}\oplus V_{1,3})\sharp\A \geq\dim \BN(Z_{1,1,1,3})\sharp\K[\Gamma]=1024$. Since $\dim\mathcal{K}=8$, the claim holds. Consequently, $\BN\cong\BN(V_{1,1}\oplus V_{1,3})$.

\end{proof}

\begin{rmk}\label{rmksplitNichols2}
Let $(i,j,k,\ell)\in\{(1,1,1,3),(2,1,2,3)\}$. Similar to Remark \ref{rmksplitNichols}, $ D_{i,j,k,\ell}\cong \BN(V_{i,j}\oplus V_{k,\ell})\sharp\A$ and $\BN(V_{i,j}\oplus V_{k,\ell})\sharp\gr\A\cong\BN(Z_{i,j,k,\ell})\sharp\K[\Gamma]$. If $i=1$ or $2$, then a direct computation shows that the Dynkin diagram of $\BN(Z_{i,j,k,\ell})$ is \Dchainthree{\xi^{-1}}{v_1}{\xi}{-1}{x}{\xi^{-1}}{\xi}{e_1} or \Dchainthree{-1}{v_1}{\xi^{-1}}{-1}{x}{\xi}{-1}{e_1}, respectively.
Moreover, $\BN(V_{i,j}\oplus V_{k,\ell})\sharp\BN(W)\cong \BN(Z_{i,j,k,\ell})$ and the Nichols algebras  in Proposition \ref{proTwoTwo} can be obtained by splitting the Nichols algebras of diagonal type.
\end{rmk}

\begin{lem}\label{lemIn}
If $V$ is isomorphic either to $V_{1,j}\oplus\K_{\chi}\oplus\K_{\chi}$, $V_{2,j}\oplus\K_{\chi^3}\oplus\K_{\chi^3}$, $V_{1,1}\oplus V_{1,3}\oplus\K_{\chi}$ or $V_{2,1}\oplus V_{2,3}\oplus\K_{\chi^3}$, then $\dim\BN(V)=\infty$.
 \end{lem}
 \begin{proof}
 Let $V_{i,j}\oplus\K_{\chi^{2i-1}}\oplus\K_{\chi^{2i-1}}=\K\{v_1,v_2\}\oplus\K\{v_3\}\oplus\K\{v_4\}\in\KYD$ for $i\in\I_{0,1},j\in\{1,3\}$. Then $V_{i,j}\oplus\K_{\chi^{2i-1}}\oplus\K_{\chi^{2i-1}}\in\AYD$. Consider the subalgebra $A_{i,j}$ of $\BN(V_{i,j}\oplus\K_{\chi^{2i-1}}\oplus\K_{\chi^{2i-1}})\sharp\A$ generated by $g,x,v_1,v_3,v_4$. A straight computation shows that $A_{i,j}$ is a pointed Hopf subalgebra with $\G(A_{i,j})\cong\Gamma$. Moreover, $\text{gr}\,A_{i,j}=\BN(T_{i,j})\sharp\K[\Gamma]$, where $T_{i,j}=\K\{x,v_1,v_3,v_4\}\in{}_{\Gamma}^{\Gamma}\mathcal{YD}$ with the Yetter-Drinfeld structure given by
 \begin{align*}
g\cdot x&=-x,\quad  g\cdot v_1=\xi^{-j}v_1,\quad g\cdot v_3=-v_3\quad g\cdot v_4=-v_4\\ \delta(x)&=g\otimes x,\quad\delta(v_1)=g^i\otimes v_1,\quad\delta(v_3)=g^{2i-1}\otimes v_3,\quad \delta(v_4)=g^{2i-1}\otimes v_4.
\end{align*}
Then the generalized Dynkin Diagram of $T_{1,j}$ or $T_{2,j}$ is given by
\begin{align*}
\xymatrix{ & \overset{-1}{\underset{v_4}{\circ}} \ar  @{-}[d]^{\xi^j} &   \\
	\overset{-1}{\underset{x}{\circ}}\ar  @{-}[r]^{\xi^j}  & \overset{\xi^j}{\underset{v_1}{\circ}} \ar  @{-}[r]^{\xi^j}  & \overset{-1}{\underset{v_3}{\circ}},\text{~or~}}
\xymatrix{ & \overset{-1}{\underset{v_4}{\circ}} \ar  @{-}[d]^{\xi^j} &   \\
	\overset{-1}{\underset{x}{\circ}}\ar  @{-}[r]^{\xi^{-j}}  & \overset{-1}{\underset{v_1}{\circ}} \ar  @{-}[r]^{\xi^j}  & \overset{-1}{\underset{v_3}{\circ}}.}
\end{align*}
Since \xymatrix{ \overset{\xi^j}{\underset{v_1}{\circ}} \ar  @{-}[r]^{\xi^j}  & \overset{-1}{\underset{v_3}{\circ}}} and \xymatrix{
	\overset{-1}{\underset{x}{\circ}}\ar  @{-}[r]^{\xi^{-j}}  & \overset{-1}{\underset{v_1}{\circ}} \ar  @{-}[r]^{\xi^j}  & \overset{-1}{\underset{v_3}{\circ}}} have infinite root systems, it follows that $\BN(V_{i,j}\oplus\K_{\chi^{2i-1}}\oplus\K_{\chi^{2i-1}})$ is infinite dimensional.

Similarly, let $V_{i,1}\oplus V_{i,3}\oplus\K_{\chi^{2i-1}}=\K\{v_1,v_2\}\oplus\K\{e_1,e_2\}\oplus\K\{v\}\in\KYD$ for $i\in\I_{0,1}$. Then $V_{i,1}\oplus V_{i,2}\oplus\K_{\chi^{2i-1}}\in\AYD$. Consider the subalgebra $A_{i}$ of $\BN(V_{i,1}\oplus V_{i,3}\oplus\K_{\chi^{2i-1}})\sharp\A$ generated by $g,x,v_1,e_1,v$. A straight computation shows that $A_{i}$ is a pointed Hopf subalgebra with $\G(A_{i})\cong\Gamma$. Moreover, $\text{gr}\,A_{i}=\BN(F_{i})\sharp\K[\Gamma]$, where $F_{i}=\K\{x,v_1,e_1,v\}\in{}_{\Gamma}^{\Gamma}\mathcal{YD}$ is a braided vector space of diagonal type. Moreover, The generalized Dynkin Diagram of $F_{1}$ or $F_{2}$ is given by
\begin{align*}
\xymatrix{ & \overset{-1}{\underset{v}{\circ}}   \ar  @{-}[dl]_{\xi}  \ar  @{-}[dr]^{\xi^{-1}}  &   \\
\overset{\xi^{-1}}{\underset{v_1}\circ}\ar  @{-}[r]^{\xi}  & \overset{-1}{\underset{x}\circ} \ar  @{-}[r]^{\xi^{-1}}  & \overset{\xi}{\underset{e_1}\circ}\text{,~or~}}
\xymatrix{ & \overset{-1}{\underset{v}{\circ}}   \ar  @{-}[dl]_{\xi}  \ar  @{-}[dr]^{\xi^{-1}}  &   \\
\overset{-1}{\underset{v_1}\circ}\ar  @{-}[r]^{\xi^{-1}}  & \overset{-1}{\underset{x}\circ} \ar  @{-}[r]^{\xi}  & \overset{-1}{\underset{e_1}\circ}.}
\end{align*}
We apply the reflection at point $v$ described as follows:
\begin{align*}
& \xymatrix{ & \overset{-1}{\underset{v}{\circ}}   \ar  @{-}[dl]_{\xi}  \ar  @{-}[dr]^{\xi^{-1}}  &   \\
\overset{\xi^{-1}}{\underset{v_1}\circ}\ar  @{-}[r]^{\xi}  & \overset{-1}{\underset{x}\circ} \ar  @{-}[r]^{\xi^{-1}}  & \overset{\xi}{\underset{e_1}\circ},}
&
&\xymatrix{\ar@/^0.5pc/@{<->}[r]^{v} & }
&
& \xymatrix{ & \overset{-1}{\underset{1 }{\circ}}   \ar  @{-}[dl]_{\xi^{-1}}  \ar  @{-}[dr]^{\xi^{}}  &   \\
\overset{-1}{\underset{2}\circ}\ar  @{-}[r]^{\xi}  & \overset{-1}{\underset{3}\circ} \ar  @{-}[r]^{\xi^{-1}}  & \overset{-1}{\underset{4 }\circ}.}
\end{align*}
Thus $\dim\BN(V_{1,1}\oplus V_{1,3}\oplus\K_{\chi^{1}})=\dim\BN(V_{2,1}\oplus V_{2,3}\oplus\K_{\chi^{3}})=\infty$ by \cite[Theorem\,17\,(p90,Step\,10)]{H09}. Indeed, one can consider subroot system $\roots(\alpha_1,\alpha_2+\alpha_3,\alpha_4)$, a direct computation shows that it is an infinite system. Consequently, we have shown that the claim holds.
\end{proof}

\begin{proofthma}

Let $V\in\KYD$ such that $\BN(V)$ is finite-dimensional. By \cite[Theorem\,4.5]{GG16}, $V$ must be semisimple.
If $V$ is the direct sum of one-dimensional objects in $\KYD$, then the theorem follows by Lemma \ref{lemNicholbyone}. Otherwise, the theorem follows by Lemma \ref{lemIn}, Propositions \ref{proNicholsTwoOfK} \& \ref{proNicholsTwoone} \& \ref{proTwoTwo}.
\end{proofthma}

\begin{rmk}
The Dynkin diagrams of rank $2$ in Remarks \ref{rmksplitNichols} \ref{rmksplitNichols1} \& \ref{rmksplitNichols2} appeared in \cite[Table\,1,\,row\,2]{H09} and the Dynkin diagrams of rank $3$ appeared in \cite[Table\,2,\,row\,8]{H09}. They are of standard type $A_2$ or $A_3$.
\end{rmk}

\begin{cor}\label{corCartan}
The generalized Cartan matrix of $\BN(M)$ is of type $A_2$, where $M$ is isomorphic either to
$V_{1,j}\oplus \K_{\chi}$, $V_{2,j}\oplus\K_{\chi^3}$, $V_{1,1}\oplus V_{1,3}$ or $V_{2,1}\oplus V_{2,3}$. Moreover, they can give rise to a Weyl-groupoid of standard type $A_2$.
\end{cor}
\begin{proof}
Note that $(\K_{\chi})^*\cong \K_{\chi^3}$ and $V_{1,j}^*\cong V_{2,j}$ by Remark \ref{rmkDmoddual}.
By Theorem \ref{thmA}, a direct computation shows that the following isomorphisms hold in $\BN(M)$, where $M$ is isomorphic either to
$V_{1,j}\oplus \K_{\chi}$, $V_{2,j}\oplus\K_{\chi^3}$, $V_{1,1}\oplus V_{1,3}$ or $V_{2,1}\oplus V_{2,3}$.
\begin{align*}
(\ad_c V_{2,j})(\K_{\chi^3})\cong V_{1,-j}\cong (\ad_c\K_{\chi^3})(V_{2,j}),\quad (\ad_c V_{1,j})(\K_{\chi})\cong V_{2,-j}\cong (\ad_c\K_{\chi})(V_{1,j}),\\
(\ad_c V_{1,1})(V_{1,3})\cong \K_{\chi^3}\cong (\ad_c V_{1,3})(V_{1,1}),\quad (\ad_c V_{2,3})(V_{2,1})\cong \K_{\chi}\cong (\ad_c V_{2,1})(V_{2,3}),\\
(\ad_c V_{2,j})^2(\K_{\chi^3})=0= (\ad_c\K_{\chi^3})^2(V_{2,j}),\quad (\ad_c V_{1,j})^2(\K_{\chi})=0= (\ad_c\K_{\chi})^2(V_{1,j}),\\
(\ad_c V_{1,1})^2(V_{1,3})=0= (\ad_c V_{1,3})^2(V_{1,1}),\quad (\ad_c V_{2,3})^2(V_{2,1})=0= (\ad_c V_{2,1})^2(V_{2,3}).
\end{align*}

Thus the generalized Cartan matrix of $\BN(M)$ is of type $A_2$, where $M$ is isomorphic either to
$V_{1,j}\oplus \K_{\chi}$, $V_{2,j}\oplus\K_{\chi^3}$, $V_{1,1}\oplus V_{1,3}$ or $V_{2,1}\oplus V_{2,3}$. Moreover, they can give rise to a Weyl-groupoid of standard type $A_2$ by \cite{AHS10,HS10,HS10b}.
\end{proof}
\begin{rmk}
We have shown that the Nichols algebras of non-diagonal type in Theorem \ref{thmA} are of type $A_1$ or standard type $A_2$ and they can be obtained by splitting the Nichols algebras $\BN(U)$ of standard type $A_2$ or $A_3$. That is, they can be characterized by restrictions of the Nichols algebras $\BN(U)$. We refer to \cite{CL17} for the details on Nichols algebras with restricted root systems and to \cite{AA} for a characterization of finite-dimensional Nichols algebras over basic Hopf algebras.
\end{rmk}
%\begin{rmk}
%We noticed that there is a recent work by N.~Andruskiewitsch and I.-E.~Angiono on Nichols algebras over basic Hopf algebras \cite{AA} and they gave a characterization of finite dimensional Nichols algebras over basic Hopf algebras $($see \cite[Theorem 1.1]{AA}$)$ while attending the workshop on Tensor categories, Hopf algebras and Quantum groups at Marburg, January 22-26, 2018.
%\end{rmk}
\section{Hopf algebras over $\mathcal{K}$}\label{secHopfalgebra}
In this section, we mainly compute the liftings of the Nichols algebras in Theorem \ref{thmA}.
We first show that the Nichols algebra $\BN(V)$ does not admit non-trivial deformations, where $V$ is isomorphic either to $\oplus_{i=1}^n\K_{\chi^{n_i}}$ with $n_i\in\{1,3\}$, $V_{2,1}$, $V_{2,3}$, $V_{2,1}\oplus\K_{\chi^3}$, $V_{2,3}\oplus\K_{\chi^3}$ or $V_{2,1}\oplus V_{2,3}$.

\begin{pro}\label{proHopfalgebranonlifting}
Let $A$ be a finite-dimensional Hopf algebra over $\mathcal{K}$ such that $\gr A\cong\BN(V)\sharp\mathcal{K}$, where $V$ is isomorphic either to $\oplus_{i=1}^n\K_{\chi^{n_i}}$ with $n_i\in\{1,3\}$, $V_{2,1}$, $V_{2,3}$, $V_{2,1}\oplus\K_{\chi^3}$, $V_{2,3}\oplus\K_{\chi^3}$ or $V_{2,1}\oplus V_{2,3}$. Then $A\cong\gr A$.
\end{pro}
\begin{proof}
We prove the proposition by showing that the defining relations of $\gr A$ also hold in $A$.

Assume that $V\cong\oplus_{i=1}^n\K_{\chi^{n_i}}=\oplus_{i=1}^n\K\{x_i\}$ with $n_i\in\{1,3\}$. Then $\BN(V)\sharp\mathcal{K}\cong\bigwedge V\sharp\mathcal{K}$. As $\Delta_A(x_i)=x_i\otimes 1+a^2\otimes x_i$, it follows that $x_i^2,x_ix_j+x_jx_i\in\Pp(A)=0$. Hence the relations in $\gr A$ also hold in $A$.

Assume that $V\cong V_{2,j}$ for $j=1,3$. Then $\BN(V_{2,j})\sharp\mathcal{K}$ is generated by the elements $x,y,a,b$ satisfying relations \eqref{eqDef1}, $ x^2=0$, $yx+\xi^jxy=0$, $y^4=0$, $ax=-xa$, $bx=-xb$, $ay+\xi ya=\xi^3xba^2$ and $by+\xi yb=xa^3$ with the coalgebra structure given by \eqref{eqDef}, $\Delta(x)=x\otimes 1+a^{-j}\otimes x+\xi(1-\xi^j)ba^{-1-j}\otimes y$ and
$\Delta(y)=y\otimes 1+a^{2-j}\otimes y-\frac{1}{2}(1+\xi^j)ba^{1-j}\otimes x$. Then a direct computation shows that
\begin{align*}
\Delta(yx+\xi^jxy)&=(yx+\xi^jxy)\otimes 1+1\otimes(yx+\xi^jxy), \\
\Delta(x^2)&=x^2\otimes 1+a^2\otimes x^2+\xi(1-\xi^j)ba\otimes (yx+\xi^jxy).
\end{align*}
It follows that the relation $yx+\xi^jxy=0$ holds in $A$ and $x^2\in\Pp_{1,a^2}(A)=\Pp_{1,a^2}(\mathcal{K})=\K\{1-a^2,ba\}$. Then $x^2=\alpha_1(1-a^2)+\alpha_2ba$ for $\alpha_1,\alpha_2\in\K$. Since $ax^2=x^2a$ and $bx^2=x^2b$, it follows that $\alpha_1=0=\alpha_2$ and the relation $x^2=0$ holds in $A$. Then $\Delta(y^4)=y^4\otimes 1+1\otimes y^4$ and whence the relation $y^4=0$ holds in $A$. Consequently, $\gr A\cong A$.

Assume that $V\cong V_{2,1}\oplus V_{2,3}$. Then $\BN(V_{2,1}\oplus V_{2,3})\sharp\mathcal{K}$ is generated by the elements $a,b,x,y,z,t$ such that $a,b,x,y$ satisfy the relations of $\BN(V_{2,1})\sharp\mathcal{K}$, $a,b,z,t$ satisfy the relations of $\BN(V_{2,3})\sharp\mathcal{K}$, and $x,y,z,t$ satisfy the relations $zx+xz=0$, $tx-zy-\xi yz+\xi xt=0$, $ty-\xi yt+\frac{1}{2}(\xi-1)xz=0$.

 Let $X=tx-zy-\xi yz+\xi xt$ and $Y=ty-\xi yt+\frac{1}{2}(\xi-1)xz$. A direct computation shows that
\begin{align}
\Delta(X)&=X\otimes 1+a^2\otimes X,\quad\Delta(Y)=Y\otimes 1+1\otimes Y+\frac{1}{2}(1-\xi)ba^3\otimes X,\label{eqV2123-1}\\
\Delta(zx+xz)&=(zx+xz)\otimes 1+1\otimes (zx+xz)+(\xi-1)ba^3\otimes X.\label{eqV2123-2}
\end{align}
Then  $X=\alpha_1(1-a^2)+\alpha_2ba$ for $\alpha_1,\alpha_2\in\K$ by \eqref{eqV2123-1}. A tedious computation on $A_{[1]}$ shows that the equation \eqref{eqV2123-2} holds if and only if $\alpha_2=0$. It follows that $zx+xz+\alpha_1(\xi-1)ba^3\in\Pp(A)=0$. Since $a(zx+xz)=(zx+xz)a$ and $ba=\xi ab$, we have that $\alpha_1=0$ and whence the relations $X=0$, $Y=0$ and $zx+xz=0$ hold in $A$. Moreover, as shown in the case $V\cong V_{2,j}$ for $j=1,3$, the relations $x^2=yx+\xi xy=y^4=0$ and $z^2=tz-\xi zt=t^4=0$ hold in $A$. Consequently, $\gr A\cong A$.

Assume that $V\cong V_{2,j}\oplus\K_{\chi^3}$ for $j=1,3$. Then $\BN(V_{2,j}\oplus\K_{\chi^3})\sharp\mathcal{K}$ is generated by $x,y,z,a,b$ such that $a,b, x, y$ satisfy the relations in $\BN(V_{2,j})\sharp\mathcal{K}$, $z,a,b$ satisfy the relations in $\BN(\K_{\chi^3})\sharp\mathcal{K}$ and $x,y,z$ satisfy $xyz+xzy+yzx+zyx=0$, $(zx)^4+(xz)^4=0$, $\frac{1}{2}\xi(1+\xi^{-j})xzx-\xi^j y^2z+(1-\xi^j)yzy+zy^2=0$. As $\Delta(x)=x\otimes 1+a^{-j}\otimes x+\xi(1-\xi^j)ba^{-1-j}\otimes y$,
$\Delta(y)=y\otimes 1+a^{2-j}\otimes y-\frac{1}{2}(1+\xi^j)ba^{1-j}\otimes x$ and $\Delta(z)=z\otimes 1+a^2\otimes z$, it follows that the relations in $\BN(V_{2,j})$ and $\BN(\K_{\chi^3})$ hold in $A$. Let $Z=\frac{1}{2}\xi(1+\xi^{-j})xzx-\xi^j y^2z+(1-\xi^j)yzy+zy^2$ and $W=xyz+xzy+yzx+zyx$, we have
\begin{align*}
\Delta(Z)=Z\otimes 1+1\otimes Z,\quad \Delta(W)=W\otimes 1+a^2\otimes W+\xi(1+\xi^j)ba\otimes Z.
\end{align*}
Then the relation $Z=0$ holds in $A$ and $W=\alpha(1-a^2)+\beta ba$ for some $\alpha,\beta\in\K$. Since $bW=Wb$ and $aW=Wa$, it follows that $\alpha=0=\beta$ and whence the relation $W=0$ holds in $A$. Moreover, $\Delta((zx)^4+(xz)^4)=((zx)^4+(xz)^4)\otimes 1+1\otimes((zx)^4+(xz)^4)$ and whence the relation $(zx)^4+(xz)^4=0$ holds in $A$. Consequently, $\gr A\cong A$.
\end{proof}

Next, we define five families of Hopf algebras $\mathfrak{A}_{1,j}(\mu)$, $\mathfrak{A}_{1,j,1}(\mu,\nu)$ for $j=1,3$ and $\mathfrak{A}_{1,1,1,3}(\mu,\nu)$ and show that they are liftings of the Nichols algebras $\BN(V_{1,j})$, $\BN(V_{1,j}\oplus\K_{\chi})$ for $j=1,3$ and $\BN(V_{1,1}\oplus V_{1,3})$, respectively.
\begin{defi}\label{defiV13}
For $j\in\{1,3\}$ and $\mu\in\K$, let $\mathfrak{A}_{1,j}(\mu)$ be the algebra generated by the elements $a,b,x,y$  satisfying the relations
\begin{gather}
a^4=1, \quad b^2=0,\quad ba=\xi ab,\quad ax=\xi xa,\quad bx=\xi xb,\quad ay+ya=\xi^{3}xba^2,\quad by+yb=xa^3,\label{eqDef1j-1}\\
x^4=0,\quad x^2+2\xi y^2=\mu(1-a^2),\quad xy+yx=\mu\xi^3ba^3.\label{eqDef1j-2}
\end{gather}
\end{defi}
$\mathfrak{A}_{1,j}(\mu)$ is a Hopf algebra with the coalgebra structure given by \eqref{eqDef} and
\begin{align}
\Delta(x)=x\otimes 1+a^{-j}\otimes x-(1+\xi^j)ba^{-1-j}\otimes y,\;
\Delta(y)=y\otimes 1+a^{2-j}\otimes y+\frac{1}{2}\xi(1-\xi^j)ba^{1-j}\otimes x.\label{eqDef1j-3}
\end{align}
\begin{rmk}
It is clear that $\gr\mathfrak{A}_{1,j}(\mu)\cong \mathfrak{A}_{1,j}(0)$ and $\mathfrak{A}_{1,j}(\mu)$ with $\mu\neq 0$ is not isomorphic as a Hopf algebra to $\mathfrak{A}_{1,j}(0)$ for $j\in\{1,3\}$.

\end{rmk}

\begin{defi}\label{defiV131}
For $j\in\{1,3\}$ and $\mu\in\K$, let $\mathfrak{A}_{1,j,1}(\mu,\nu)$ be the algebra generated by the elements $a,b,x,y,z$ such that  $a,b,x,y$ satisfy \eqref{eqDef1j-1} \eqref{eqDef1j-2}, and all together they satisfy
\begin{gather*}
z^2=0,\quad az=\xi za,\quad bz=\xi zb,\quad zx^2+(1-\xi^j)xzx-\xi^j x^2z=2\nu ba^3,\\
 \xi^j xyz-\xi^j xzy+yzx+zxy=\nu(1-a^2),\quad \frac{1}{2}\xi(1+\xi^{-j})(xz)^2(yz)^2+(yz)^4+(zy)^4=\alpha(\nu).
\end{gather*}
\end{defi}
$\mathfrak{A}_{1,j,1}(\mu)$ is a Hopf algebra with the coalgebra structure given by \eqref{eqDef} \eqref{eqDef1j-3} and $\Delta(z)=z\otimes 1+a^2\otimes z$.
\begin{rmk}
\begin{itemize}
  \item Due to complicated commutators relations and computations, the explicit form of $\alpha(\nu)$ is not clear. However, if $\nu=0$, then a tedious computation shows that $\alpha(\nu)=0$ and we shall show that $\mathfrak{A}_{1,j,1}(\mu,0)$ is a non-trivial lifting of $\BN(V_{1,j}\oplus\K_{\chi})$. Moreover, assume that $\nu\neq 0$, then we shall show that $\alpha(\nu)\neq 0$ if $\mathfrak{A}_{1,j,1}(\mu,\nu)$ is a non-trivial lifting. Hence the liftings of $\BN(V_{1,j}\oplus\K_{\chi})$ are characterized by the parameters $\mu,\nu$.
  \item It is clear that $\gr\mathfrak{A}_{1,j,1}(\mu,\nu)\cong\mathfrak{A}_{1,j,1}(0,0)$ and $\mathfrak{A}_{1,j}(\mu)$ is a Hopf subalgebra of $\mathfrak{A}_{1,j,1}(\mu,0)$.
\end{itemize}

%\begin{align*}
%zyzx-\xi^jzxzy-\xi^jyzxz-xzyz=\nu(1-\xi^j)z(1-\xi^ja^2).\\
%(zx)^2+\xi^j(xz)^2=2\xi^j\nu zba^3
%\end{align*}
\end{rmk}

\begin{defi}\label{defiV1313}
For $\mu,\nu\in\K$, let $\mathfrak{A}_{1,1,1,3}(\mu,\nu)$ be the algebra generated by the elements $a,b,x,y,z,t$ such that $a,b,x,y$ satisfy \eqref{eqDef1j-1} \eqref{eqDef1j-2}, and all together they satisfy.
\begin{gather*}
az=\xi za,\quad bz=\xi zb,\quad at+ta=\xi^{3}xba^2,\quad bt+tb=xa^3,\quad z^4=0,\quad
z^2+2\xi t^2=\nu(1-a^2),\\ zt+tz=\nu\xi^3ba^3,\quad
zx-\xi xz=0,\quad tx+zy+yz+xt=0,\quad ty+yt+\frac{1}{2}(1-\xi)xz=0.
\end{gather*}
\end{defi}
$\mathfrak{A}_{1,1,1,3}(\mu,\nu)$ is a Hopf algebra with the coalgebra structure given by \eqref{eqDef} and
\begin{align}
\Delta(x)=x\otimes 1+a^{3}\otimes x-(1+\xi)ba^{2}\otimes y,\;
\Delta(y)=y\otimes 1+a \otimes y+\frac{1}{2}(1+\xi )b  \otimes x,\label{eqV1113-1}\\
\Delta(z)=z\otimes 1+a \otimes z+(\xi-1)b  \otimes t,\;
\Delta(t)=t\otimes 1+a^{3}\otimes t+\frac{1}{2}(\xi-1)ba^{2}\otimes z.\label{eqV1113-2}
\end{align}

\begin{rmk}
It is clear that $\gr\mathfrak{A}_{1,1,1,3}(\mu,\nu)\cong \mathfrak{A}_{1,1,1,3}(0,0)$ and $\mathfrak{A}_{1,1}(\mu)$, $\mathfrak{A}_{1,3}(\nu)$ are Hopf subalgebras of $\mathfrak{A}_{1,1,1,3}(\mu,\nu)$.
\end{rmk}

\begin{lem}\label{lemPBW}
Let $j\in\{1,3\}$. A linear basis of $\mathfrak{A}_{1,j}(\mu)$, $\mathfrak{A}_{1,j,1}(\mu,0)$ or $\mathfrak{A}_{1,1,1,3}(\mu,\nu)$ is given by $\{y^ix^jb^ka^{\ell},i,k\in\I_{0,1},j,\ell\in\I_{0,3}\}$, $\{z^{n_1}(yz)^{n_2}(xz)^{n_3}y^{n_4}x^{n_5}b^{n_6}a^{n_7},\;\;n_2,n_5,n_7\in\I_{0,3},n_1,n_3,n_4,n_6\in\I_{0,1}\}$ or\\ $\{t^{n_1}z^{n_2}(yt)^{n_{3}}y^{n_4}x^{n_5}b^{n_6}a^{n_7}:\;n_2,n_4,n_7\in\I_{0,3},n_1,n_{3},n_5,n_6\in\I_{0,1}\}$, respectively.
\end{lem}
\begin{proof}
We first prove the statement for $\mathfrak{A}_{1,j}(\mu)$ by applying the Diamond Lemma \cite{B} with the order $y<x<b<a$. By the Diamond Lemma, it suffices to show that all overlap ambiguities are resolvable. That is, the ambiguities can be reduced to the same expression by different substitution rules. To verify all the ambiguities are resolvable is tedious but straightforward. Here we only check the overlaps $(xy)y=x(y^2)$, $x^3(xy)=(x^4)y$ and $(ay)y=a(y^2)$.
 Note that $ax^2=-x^2a$. After a direct computation, $ba^3y-yba^3=-xa^2$. Then
\begin{align*}
(xy)y&=-yxy+\mu\xi^3ba^3y=-y(-yx+\mu\xi^3ba^3)+\mu\xi^3ba^3y\\
&=y^2x+\mu\xi^3(ba^3y-yba^3)=y^2x+\mu\xi xa^2\\
&=\frac{1}{2}\xi(x^2-\mu(1-a^2))x+\mu\xi xa^2=\frac{1}{2}\xi x^3-\frac{1}{2}\xi\mu x+\frac{1}{2}\xi\mu xa^2
=\frac{1}{2}\xi x^3+\frac{1}{2}\xi^3\mu x(1-a^2)=x(y^2).
\end{align*}
Note that $ba^3x=xba^3$. Then $x(xy)=yx^2$ and whence $x^3(xy)=yx^4=0=(x^4)y$. Similarly, we have that
\begin{align*}
(ay)y&=(-ya+\xi^3xba^2)y=-yay+\xi^3xba^2y
=-y(-ya+\xi^3xba^2)+\xi^3x(-yba^2+xa)\\
&=y^2a+\xi(xy+yx)ba^2+\xi^3x^2a=y^2a+\xi^3x^2a\\
&=\frac{1}{2}\xi x^2a+\frac{1}{2}\mu\xi^3(1-a^2)a+\xi^3x^2a=\frac{1}{2}\xi ax^2+\frac{1}{2}\xi^3\mu a(1-a^2)=a(y^2).
\end{align*}

For $\mathfrak{A}_{1,j,1}(\mu,0)$ and $\mathfrak{A}_{1,1,1,3}(\mu,\nu)$, the proof follows the same line as for $\mathfrak{A}_{1,j}(\mu)$ and we omit the tedious details.

Assume that $\nu\neq 0$ for $\mathfrak{A}_{1,j,1}(\mu,\nu)$, a tedious computation shows that $\alpha(\nu)\neq 0$, otherwise the ambiguities like $(a(yz))(yz)^3=a(yz)^4$ are not resolvable.
\end{proof}

\begin{lem}
For $j=1,3$, $\gr\mathfrak{A}_{1,j}(\mu)\cong\BN(V_{1,j})\sharp\mathcal{K}$, $\gr\mathfrak{A}_{1,j,1}(\mu,0)\cong\BN(V_{1,j}\oplus\K_{\chi})\sharp\mathcal{K}$ and $\gr\mathfrak{A}_{1,1,1,3}(\mu)\cong\BN(V_{1,1}\oplus V_{1,3})\sharp\mathcal{K}$.
\end{lem}
\begin{proof}
Let $\Lambda_0$ be the subalgebra of $\mathfrak{A}_{1,j}(\mu)$ generated by the elements $a,b$. Note that $\dim\Lambda_0=8$ by Lemma \ref{lemPBW}. Then it is easy to see that $\Lambda_0\cong\mathcal{K}$ as Hopf algebras. Let $\Lambda_1=\Lambda_0+\mathcal{K}\{x,y\}$,  $\Lambda_2=\Lambda_1+\mathcal{K}\{x^2,xy\}$, $\Lambda_3=\Lambda_2+\mathcal{K}\{x^3,x^2y\}$ and $\Lambda_4=\Lambda_3+\mathcal{K}\{x^3y\}$. A direct computation shows that $\{\Lambda_{\ell}\}_{\ell=0}^4$ is a coalgebra filtration of $\mathfrak{A}_{1,j}(\mu)$. Hence $(\mathfrak{A}_{1,j}(\mu))_0\subseteq \mathcal{K}$ and $(\mathfrak{A}_{1,j}(\mu))_{[0]}\cong \mathcal{K}$, i.e., $\mathfrak{A}_{1,j}(\mu)$ is a Hopf algebra over $\mathcal{K}$. Hence, $\gr\mathfrak{A}_{1,j}(\mu) \cong R_{1,j}\sharp H$. It is clear that $\Pp(R_{1,j})\cong V_{1,j}$ in $\KYD$ by the definition of $\mathfrak{A}_{1,j}(\mu)$. Moreover, $\dim R_{1,j}=8=\dim\BN(V_{1,j})$ by Lemma \ref{lemPBW}. It follows that $\gr\mathfrak{A}_{1,j}(\mu)\cong\BN(V_{1,j})\sharp\mathcal{K}$.

For $\mathfrak{A}_{1,j,1}(\mu,0)$ and $\mathfrak{A}_{1,1,1,3}(\mu,\nu)$, the proof follows the same line as for $\mathfrak{A}_{1,j}(\mu)$.
\end{proof}

\begin{pro}\label{proHopfalgebraliftings}
Let $A$ be a finite-dimensional Hopf algebra over $\mathcal{K}$ such that $\gr A\cong\BN(V)\sharp\mathcal{K}$, where $V$ is isomorphic either to $V_{1,j}$ for $j=1,3$,   $V_{1,j}\oplus\K_{\chi}$ for $j=1,3$, or $V_{1,1}\oplus V_{1,3}$. Then $A$ is isomorphic either to
$\mathfrak{A}_{1,j}(\mu)$, $\mathfrak{A}_{1,j,1}(\mu,\nu)$ or $\mathfrak{A}_{1,1,1,3}(\mu,\nu)$ for $\mu,\nu\in\K$.
\end{pro}
\begin{proof}
Assume that $V\cong V_{1,j}$ for $j=1,3$. Note that $\gr A\cong\mathfrak{A}_{1,j}(0)$. Then by  \eqref{eqDef1j-3}, a direct computation shows that
\begin{gather}
%\Delta(x^4)=x^4\otimes 1+1\otimes x^4,\;
\Delta(x^2+2\xi y^2)=(x^2+2\xi y^2)\otimes 1+a^2\otimes (x^2+2\xi y^2),\label{eqV1jC-1}\\
\Delta(xy+yx)=(xy+yx)\otimes 1+1\otimes (xy+yx)+\xi ba^3\otimes (x^2+2\xi y^2).\label{eqV1jC-2}
\end{gather}
By \eqref{eqV1jC-1}, $x^2+2\xi y^2\in\Pp_{1,a^2}(A)=\Pp_{1,a^2}(\mathcal{K})=\K\{1-a^2,ba\}$. It follows that $x^2+2\xi y^2=\alpha_1(1-a^2)+\alpha_2ba$ for $\alpha_1,\alpha_2\in\K$. Then by \eqref{eqV1jC-2}, $\Delta(xy+yx+\alpha_1\xi ba^3)=(xy+yx+\alpha_1\xi ba^3)\otimes 1+1\otimes (xy+yx+\alpha_1\xi ba^3)+\alpha_2\xi ba^3\otimes ba$. A direct computation on $A_{[1]}$ shows that the equation \eqref{eqV1jC-2} holds if and only if $\alpha_2=0$. It follows that $x^2+2\xi y^2=\alpha_1(1-a^2)$ and $xy+yx+\alpha_1\xi ba^3=0$ hold in $A$. Moreover, $\Delta(x^4)=x^4\otimes 1+1\otimes x^4$. Thus there is an epimorphism from $\mathfrak{A}_{1,j}(\alpha_1)$ to $A$. Since $\dim A=\dim \mathfrak{A}_{1,j}(\alpha_1)$ by Lemma \ref{lemPBW},  $A\cong\mathfrak{A}_{1,j}(\alpha_1)$ .

Assume that $V\cong V_{1,j}\oplus\K_{\chi}$ for $j=1,3$. Note that $\gr A\cong\mathfrak{A}_{1,j,1}(0,0)$. As shown in the case $V\cong V_{1,j}$, $x^4=z^2=0$, $x^2+2\xi y^2=\mu(1-a^2)$, $xy+yx=\mu\xi^3ba^3$ for some $\mu\in\K$.
Let $L=zx^2+(1-\xi^j)xzx-\xi^j x^2z$  and $M=\xi^j xyz-\xi^j xzy+yzx+zxy$.
By  \eqref{eqDef1j-3} and $\Delta(z)=z\otimes 1+a^2\otimes z$, a direct computation shows that
\begin{align*}
\Delta(z^2)=z^2\otimes 1+1\otimes z^2,\quad\Delta(M)=M\otimes 1+a^2\otimes M,\quad\Delta(L)=L\otimes 1+1\otimes L-2ba^3\otimes M.
\end{align*}
Then as stated above, $z^2=0$, $M=\nu(1-a^2)$ and $L=2\nu ba^3$ for some $\nu\in\K$. Let $N= \frac{1}{2}\xi(1+\xi^{-j})(xz)^2(yz)^2+(yz)^4+(zy)^4$. Then $\Delta(N)=N\otimes 1+1\otimes N+\alpha$, where $\alpha\in \sum_{i=1}^{7}A_{[i]}\otimes A_{[8-i]}$. Assume that $\nu=0$, then a tedious computation shows that $\Delta(N)=N\otimes 1+1\otimes N$ and whence the relation $N=0$ holds in $A$. Thus by Lemma \ref{lemPBW}, the liftings of $\BN(V)$ have the form $\mathfrak{A}_{1,j,1}(\mu,\nu)$.

Assume that $V\cong V_{1,1}\oplus V_{1,3}$. Note that $\gr A\cong\mathfrak{A}_{1,1,1,3}(0)$.
Let $X=tx+zy+yz+xt$ and $Y=ty+yt+\frac{1}{2}(1-\xi)xz$. By \eqref{eqV1113-1} \eqref{eqV1113-2}, we have
\begin{align}
\Delta(X)=X\otimes 1+a^2\otimes X,\quad \Delta(Y)=Y\otimes 1+1\otimes Y+\frac{1}{2}(\xi-1)ba^3\otimes X,\label{eqCV1113-1}\\
\Delta(zx-\xi xz)=(zx-\xi xz)\otimes 1+1\otimes(zx-\xi xz)+(\xi-1)ba^3\otimes X. \label{eqCV1113-2}
\end{align}
Then  $X=\alpha_1(1-a^2)+\alpha_2ba$ for $\alpha_1,\alpha_2\in\K$ by \eqref{eqCV1113-1}. A tedious computation on $A_{[1]}$ shows that the equation \eqref{eqCV1113-2} holds if and only if $\alpha_2=0$. It follows that $zx-\xi xz+\alpha_1(\xi-1)ba^3\in\Pp(A)=0$. Since $a(zx-\xi xz)=-(zx-\xi xz)a$ and $ba=\xi ab$, we have that $\alpha_1=0$ and whence the relations $X=0$, $Y=0$ and $zx-\xi xz=0$ hold in $A$. Moreover, as shown in the case $V\cong V_{1,j}$, $x^4=z^4=0$, $x^2+2\xi y^2=\mu(1-a^2)$, $xy+yx=\mu\xi^3ba^3$, $z^2+2\xi t^2=\nu(1-a^2)$ and $xy+yx=\nu\xi^3ba^3$ for some $\mu,\nu\in\K$.  Thus there is an epimorphism from $\mathfrak{A}_{1,1,1,3}(\mu,\nu)$ to $A$. Since $\dim A=\dim \mathfrak{A}_{1,1,1,3}(\mu,\nu)$ by Lemma \ref{lemPBW}, $A\cong \mathfrak{A}_{1,1,1,3}(\mu,\nu)$.
\end{proof}

\begin{thm}\label{thmGeneratedbyone}
Let $A$ be a finite-dimensional Hopf algebra over $\mathcal{K}$ with the infinitesimal braiding $V$. Assume that the diagram of $A$ is strictly graded. Then $\gr A\cong\BN(V)\sharp \mathcal{K}$.
\end{thm}
\begin{proof}
Since $A$ is a Hopf algebra over $\mathcal{K}$, $\gr A\cong R\sharp\mathcal{K}$, where $R$ is the diagram of $A$. Let $S$ be the graded dual of $R$. Then $S$ is generated by $S(1)$ if and only if $\Pp(R)=R(1)$ by \cite[Lemma\,2.4]{AS02}. Since $R$ is strictly graded, there is an epimorphism $S\twoheadrightarrow \BN(W)$, where $W=S(1)$. To show that $R$ is generated by $R(1)$,
it suffices to show that the defining relations of $\BN(W)$ hold in $S$.

Assume $W=\oplus_{k=1}^n\K_{\chi^{i_k}}:=\oplus_{k=1}^n\K\{v_k\}$ with $i_k\in\{1,3\}$. By Lemma \ref{lemNicholbyone}, $\BN(W)=\bigwedge W$. Then it is easy to see that $r$ is a primitive element and $c(r\otimes r)=r\otimes r$  for any defining relation $r$ of $\BN(W)$,  since $c=-\tau$, where $\tau$ is a flip.

Assume that $W=V_{1,j}$ with $j\in\{1,3\}$. By Proposition \ref{proNicholsTwoOfK},
$\BN(W)$ is generated by  $v_1, v_2$ satisfying the relations $v_1^4=0$, $v_1^2+2\xi v_2^2=0$ and $v_1v_2+v_2v_1=0$ and these defining relations are primitive in $S$. Hence it suffices to show that $c(r\otimes r)=r\otimes r$ for $r=v_1^4$, $v_1^2+2\xi v_2^2$ and $v_1v_2+v_2v_1$. As
$\delta(v_1)=a^{-j}\otimes v_1-(1+\xi^j)ba^{-1-j}\otimes v_2$ and
$\delta(v_2)=a^{2-j}\otimes v_2+\frac{1}{2}\xi(1-\xi^j)ba^{1-j}\otimes v_1$, we have that
\begin{gather*}
\delta(v_1^4)=1\otimes v_1^4,\;
\delta(v_1^2+2\xi v_2^2)=a^2\otimes (v_1^2+2\xi v_2^2),
\delta(v_1v_2+v_2v_1)=1\otimes (v_1v_2+v_2v_1)+\xi ba^3\otimes (v_1^2+2\xi v_2^2).
\end{gather*}
Then it is easy to see that the claim follows. Similarly, the claim follows for the remaining cases.
\end{proof}

Finally, we give a proof of Theorem \ref{thmB}.
\begin{proofthmb}
These Hopf algebras appearing in the theorem are pairwise non-isomorphic since their infinitesimal braidings are pairwise non-isomorphic as Yetter-Drinfeld modules over $\mathcal{K}$. Others follow by Propositions \ref{proHopfalgebranonlifting} \& \ref{proHopfalgebraliftings} and Theorem \ref{thmGeneratedbyone}.
\end{proofthmb}

\vskip10pt \centerline{\bf ACKNOWLEDGMENT}

\vskip10pt This paper was written during the visit of the author to University of Padova supported by China Scholarship Council (No.~201706140160) and the NSFC (Grant No.~11771142). The author is indebted to his supervisors Profs. Naihong Hu and Giovanna Carnovale for their kind help and continued encouragement. The author is grateful to Prof. G.-A. Garcia for helpful conversations and his invaluable encouragement. The author would like to thank Profs. G.-A. Garcia and L. Vendramin for providing the computation with GAP of the Nichols algebras, and thank Profs. Andruskiewitsch and Angiono very much for helpful comments and discussion.

%\end{CJK}
\end{document}